\documentclass{article}

\usepackage{tikz}
\usepackage{mathrsfs}
\usepackage{amssymb}
\usepackage{epsfig}
\usepackage{amsmath}

\newtheorem{Theorem}{Theorem}
\newtheorem{Lemma}{Lemma}
\newtheorem{Proposition}[Lemma]{Proposition}
\newtheorem{Corollary}[Lemma]{Corollary}
\newtheorem{Definition}{Definition}

\newcommand{\spg}{\mathrm{Spin}^{G}}
\newcommand{\tf}{\widetilde{F}^{+k}(M)}
\newcommand{\tfx}{\widetilde{F}_{x}^{+k}(M)}

\newcommand{\tsln}{\widetilde{\mathrm{SL}}(\mathbf{R}^{n})}

\newcommand{\gln}{\mathrm{GL}(\mathbf{R}^{n})}
\newcommand{\glnp}{\mathrm{GL}^{+}(\mathbf{R}^{n})}
\newcommand{\glnvp}{\mathrm{GL}^{+}(\mathbf{R}^{4})}
\newcommand{\tglnvp}{\widetilde{\mathrm{GL}}{}^{+}(\mathbf{R}^{4})}
\newcommand{\tglnp}{\widetilde{\mathrm{GL}}{}^{+}(\mathbf{R}^{n})}
\newcommand{\lsln}{\mathfrak{sl}(\mathbf{R}^{n})}
\newcommand{\id}{\mathrm{id}}

\newcommand{\dg}{\mathrm{dim}(\mathfrak{g})}
\newcommand{\mdiv}{\mathrm{Div}}
\newcommand{\proof}{{\bf Proof \,}}

\newcommand{\del}{\partial}
\newcommand{\supp}{\mathrm{Supp}}
\newcommand{\qed}{${}$ $\hfill \square$}
\newcommand{\vecn}{\mathrm{Vec}_{n}}
\newcommand{\one}{\mathbf{1}}
\begin{document}
\title{Bundles with a lift of infinitesimal diffeomorphisms}
\author{Bas Janssens}
\date{\small{University of Utrecht}}
\maketitle
\begin{center}
\end{center}
\begin{abstract}
\noindent We slightly extend the notion of a natural fibre bundle by requiring
diffeomorphisms of the base to lift to automorphisms of the bundle
only infinitesimally, i.e. at the level of the Lie algebra of vector fields.
Spin structures are natural only in this extended sense.
We classify fibre bundles with this property, assuming a 
finite dimensional structure group.
This includes all spin structures, but only some $\mathrm{Spin}^{c}$
and $\spg$-structures. This classification links the gauge group $G$
to the topology of space-time.
\end{abstract}

\section{Introduction}\label{een}
All bundles are equal, but some are
more natural than others. 
For example, the naturality of
tangent bundles,
frame bundles and jet bundles is beyond all question.
A bundle is called `natural' if diffeomorphisms of the base 
lift to automorphisms of the bundle in a local fashion. 

Let us phrase this more carefully.
Any principal fibre bundle $\pi : P \rightarrow M$ determines a 
sequence of
groups
\begin{eqnarray}\label{flapperdrop}
1 \rightarrow 
\Gamma_{c}(\mathrm{Ad}(P))
\rightarrow
\mathrm{Aut}_{c}(P)
\rightarrow
\mathrm{Diff}_{c}(M) \rightarrow 1 \,, 
\end{eqnarray}
where $\mathrm{Diff}_{c}(M)$ is the group of 
compactly
supported diffeomorphisms of $M$, 
$\mathrm{Aut}_{c}(P)$ is its preimage 
in the group of automorphisms of $P$, and we have identified the gauge 
group of vertical automorphisms with sections of the adjoint bundle 
$\mathrm{Ad}(P)$. 
 
\begin{Definition}
A `natural' principal fibre bundle is
a principal fibre bundle $P$ for which (\ref{flapperdrop})
is exact and split, together with 
a distinguished splitting 
homomorphism
$\Sigma : \mathrm{Diff}_{c}(M) \rightarrow \mathrm{Aut}_{c}(P)$.
Moreover, 
$\Sigma$ is required to be local in the sense that 
it should lift each diffeomorphism 
$\phi : U_1 \rightarrow U_2$
between open subsets of $M$
to a bundle isomorphism 
$\Sigma(\phi) : \pi^{-1}(U_1) \rightarrow \pi^{-1}(U_2)$
in a functorial fashion.
\end{Definition}

Natural bundles have been classified.
A theorem of Salvioli, Palais, Terng, Ep\-stein and Thurston
(see \cite{Sa}, \cite{PT} and \cite{ET}) states that 
any natural fibre bundle is associated to 
the $k^{\mathrm{th}}$ order frame bundle $F^{k}(M)$.

In this paper, we seek to extend the notion of a natural fibre bundle
in two separate ways. First of all, we do not require locality, but prove it.
And secondly, we only require diffeomorphisms of the base 
to lift to automorphisms of the bundle infinitesimally, i.e.
at the level of Lie algebras. We will call the principal fibre bundles 
which are natural in this extended sense `infinitesimally natural',
and our main objective 
will be their classification.
 
Let us formulate this more precisely. For any principal fibre bundle, 
the sequence of groups (\ref{flapperdrop}) gives rise 
to the exact sequence of Lie algebras
\begin{equation}\label{rijtje}
0 \rightarrow \Gamma_{c}(\mathrm{ad}(P)) \rightarrow \Gamma_{c}(TP)^{G}
\rightarrow \Gamma_{c}(TM) \rightarrow 0\,,
\end{equation}
where `$c$' stands for `$0$ outside a compact subset of $M$'.

\begin{Definition}
An `infinitesimally natural' principal fibre bundle is
a principal fibre bundle $P$, together with 
a distinguished 
Lie algebra homomorphism  
$\sigma :\Gamma_{c}(TM) \rightarrow \Gamma_{c}(TP)^{G}$
that splits the exact sequence of Lie algebras (\ref{rijtje}).
\end{Definition}

We emphasise that we will not require $\sigma$ to be local, 
continuous or $C^{\infty}(M)$-linear, and it need not come from a map
of bundles. 
We only require $\sigma$ to be a 
homomorphism of Lie algebras.

The price we pay for this level of generality is that we must 
confine ourselves to fibre bundles with 
a finite dimensional structure group. This allows us to
first study principal fibre bundles, and then generalise 
the results to more general fibre bundles with a finite dimensional structure 
group, such as vector bundles.

There are three main reasons for wanting to extend the notion
of a natural bundle. First of all, natural bundles are meant
to describe geometric objects, but not all geometry is local.
One might feel that the locality requirement on $\Sigma$ is therefore
inappropriate.
The universal cover $\tilde{M} \rightarrow M$ for instance is
not a natural bundle, although 
its ties with the global topology of $M$ are unmistakable.
It does however allow for a lift of infinitesimal
diffeomorphisms, and is therefore infinitesimally natural. 
 
The second reason comes from field theory. 
Fields are described by sections 
of a fibre bundle $E$ over space-time $M$, where $E$ is 
associated to some principal fibre bundle $P$. 
Fermions, for example, 
are sections of a spinor bundle $E$, 
associated to a spin structure $P$.
As spinors transform according to a projective representation of the 
Lorentz group, they acquire a minus sign upon a full rotation.
This means that there is no hope of 
lifting global diffeomorphisms, so that spinor bundles cannot be 
natural bundles.
One might feel however that
it would be unfair to discard them as unnatural.
After all, they do occur in nature, or at least in our 
mathematical description of it. 
As one can lift diffeomorphisms infinitesimally, i.e. 
at the level of vector fields, a spinor bundle does constitute
an infinitesimally natural bundle.

The third reason is that in field theory, 
one would like to construct a 
stress-energy-momentum tensor. 
This is the conserved quantity corresponding to infinitesimal
transformations of space-time $M$.
Unfortunately, $\mathrm{Diff}(M)$
does not act on the fields $\Gamma(E)$
directly. Since $\mathrm{Aut}(P)$ does, we
need an infinitesimal lift of  
$\mathrm{Diff}(M)$ into $\mathrm{Aut}(P)$,
i.e. a Lie algebra homomorphism $\sigma$
that splits (\ref{rijtje}).
This means that an infinitesimally natural bundle is precisely 
what one needs in order to construct a SEM-tensor from Noether's 
principle. (See e.g. \cite{GM}.)
The splitting $\sigma$ is to be interpreted as the transformation 
behaviour of the fields under
infinitesimal space-time transformations.

The outline of the paper is as follows.
Sections \ref{twee}, \ref{drie} and \ref{vier} are devoted to 
the classification of infinitesimally natural principal fibre bundles.
The central result is theorem \ref{grondsmakreefknol}, which 
states the following.
\begin{center}
\emph{Any infinitesimally natural principal fibre bundle
is associated to the universal cover of the connected component of the 
$k^{\mathrm{th}}$ order frame bundle $\tf$.}
\end{center}
We extend this to fibre bundles
with a finite dimensional structure group in 
section \ref{zes}, with special attention for vector
bundles.
In section \ref{vijf}, 
we seek conditions under which 
a splitting of (\ref{rijtje}) gives rise to a flat connection.

Finally, in section \ref{zeven}, we study
spin structures in the presence of a gauge field.
Ordinary spin structures are infinitesimally natural, 
as they carry a canonical 
splitting of (\ref{rijtje}).
In the presence of a gauge group $G$ however, one should consider 
$\spg$-structures instead (see \cite{HP}, \cite{AI}). 
For example, if $G = U(1)$, the appropriate bundles are 
$\mathrm{Spin}^{c}$-structures. 
In contrast to ordinary spin structures, not all 
$\spg$-structures are infinitesimally natural.
We show that for a compact gauge group, 
infinitesimally natural $\spg$-structures 
correspond precisely to
homomorphisms
$
\pi_{1}(F(M)) \rightarrow G
$ that are injective on $\pi_{1}(\gln)$.

In a sense, this shows that complicated manifolds
call for complicated gauge groups.
Some manifolds, such as $CP^{2}$, do not allow for any 
infinitesimally natural $\spg$-structure at all.  
If one feels that fermions should have a well defined 
transformation behaviour under infinitesimal space-time
transformations, then 
these manifolds are disqualified as models for space-time.

\section{Principal bundles as Lie algebra extensions}\label{twee}

We seek to classify infinitesimally natural principal fibre bundles.
It is to this end that we study Lie algebra homomorphisms
$\sigma$ that split (\ref{rijtje}).
In this section, we will prove that $\sigma$ must be a differential
operator of finite order.

The first step, to be taken in section \ref{ideaal}, is to show 
that maximal ideals in $\Gamma_{c}(TM)$ correspond precisely 
to points in $M$. Using this, we will see in section \ref{localmap} that 
$\sigma$ must be a local map. 
Once again, we stress that locality of $\sigma$ will be a theorem, 
not an assumption.
We will then prove, in section \ref{diffop}, that $\sigma$ is in fact a
differential operator.

Before we proceed, let us have a closer look at
the exact sequence of Lie algebras 
(\ref{rijtje}), derived 
from a smooth principal $G$-bundle $\pi : P \rightarrow M$.
The last term, $\Gamma_{c}(TM)$, is the Lie algebra of smooth, compactly supported 
vector fields on $M$.
The middle term, $\Gamma_{c}(TP)^{G}$, is the Lie algebra of 
$G$-invariant vector fields $v$ 
on $P$ such that $\pi(\mathrm{Supp}(v))$ is compact.
The pushforward defines a
Lie algebra homomorphism $\pi_{*} : \Gamma_{c}(TP)^{G} \rightarrow
\Gamma_{c}(TM)$ because of $G$-invariance, and
its kernel $\Gamma_{c}(TP)_{v}^{G}$ is the ideal of vertical
vector fields. 
It can be identified with $\Gamma_{c}(\mathrm{ad}(P))$, 
the Lie algebra of compactly supported
sections of the adjoint bundle $\mathrm{ad}(P) := P \times_{G}\mathfrak{g}$,
where $\mathfrak{g}$ is the Lie algebra of $G$.


\subsection{Ideals of the Lie algebra of vector fields}\label{ideaal}
The following lemma, due to Shanks 
and Pursell \cite{SP},
constitutes the linchpin of the proof. 
It identifies the maximal ideals of $\Gamma_{c}(TM)$
as points in $M$.
The proof is taken from
\cite{SP}, with perhaps some minor clarifications. 
\begin{Lemma} \label{Ranzig}
Let $\Gamma_{c}(TM)$ be the Lie algebra
of smooth compactly supported vector fields on $M$. 
Then the maximal ideals of $\Gamma_{c}(TM)$
are labelled by points $q$ in $M$, 
$I_{q}$ being the ideal of vector fields in $\Gamma_{c}(TM)$ which are 
zero and flat at $q$.
That is,  
$I_{q} = \{ v \in \Gamma_{c}(TM) \,|\, v(q)=0 \,\, \mathrm{and}\,\, 
(\mathrm{ad}(w_{i_1})\ldots \mathrm{ad}(w_{i_n}) v)(q)=0
\, \,\forall \,\,w_{i_1} , \ldots w_{i_n} ,\in \Gamma_{c}(TM)\}\,.$

\end{Lemma}
\proof  
Let $I \subset \Gamma_{c}(TM)$ be an ideal. Suppose that 
there exists a point $q \in M$ such that $v(q) = 0$ for all $v \in
I$. Then $I \subseteq I_{q}$.
Indeed, for all $w_{i_1} , \ldots , w_{i_n} \in \Gamma_{c}(TM)$,
one has $(\mathrm{ad}(w_{i_1})\ldots \mathrm{ad}(w_{i_n}) v)(q)=0$ 
because the l.h.s. is in $I$.

Now suppose that $I$ does not have such a point. 
We will prove that this implies $I = \Gamma_{c}(TM)$.
From this, the lemma will follow. If $I$ is a maximal ideal, 
then surely $I \neq \Gamma_{c}(TM)$, so we must have $I \subseteq I_{q}$.
But it is clear from its definition that $I_{q}$ is an ideal, so, 
$I$ being maximal, we must have $I = I_{q}$. 
Conversely, the ideals
$I_{q}$ are all maximal, because any enveloping ideal must either be
all of $\Gamma_{c}(TM)$, or else be  
contained some $I_{\tilde{q}}$, which cannot be unless $\tilde{q} = q$.  

We are therefore left to prove the following: suppose that
$I$ is an ideal such that 
for all $q \in M$, there exists a $v \in I$ with $v(q) \neq 0$.
Then $I = \Gamma_{c}(TM)$. 

We will shortly prove the following statement:
for all $q \in M$, there exists a neighbourhood
$U_q $ such that all $w \in \Gamma_{c}(TM)$
with $\mathrm{Supp}(w) \subset U_q$ can be written as $w = [v,u]$
with $v \in I$ and $u \in \Gamma_{c}(TM)$.

Such $w$ are therefore elements of $I$. 
Membership of $I$ then extends to arbitrary 
$w \in \Gamma_{c}(TM)$ by covering its support with a finite 
amount of neighbourhoods 
$U_{q_1} \ldots U_{q_{N}}$, for which we construct a partition of unity
$\sum_{s=1}^{N} f_s = 1$, $\supp (f_s) \subseteq U_{q_s}$.
We may then write $w = \sum_{s=1}^{N} f_s w = \sum_{s=1}^{N} [v_s , u_s] \in I$.
This concludes the proof modulo lemma \ref{techniek}, which is next in line. \qed


Although one can easily find $u \in \Gamma_{c}(TM)$ and $v \in I$ such that 
$w = [u,v]$ locally, 
it is not always
clear how to extend these to global vector fields while simultaneously
satisfying $w = [u,v]$. For instance, with $M$ the circle 
$S^1$ and $v = w = \del_{\theta}$, the solution $u = \theta \del_{\theta}$
does not globally exist. 
We need to do some work in order to define a proper cut-off procedure.

\begin{Lemma}\label{techniek}
Let $I$ be a maximal ideal of $\Gamma_{c}(TM)$ containing a 
vector field which does not vanish at $q \in M$.
Then $q$ has a neighbourhood
$U_q $ such that all $w \in \Gamma_{c}(TM)$
with $\mathrm{Supp}(w) \subset U_q$ can be written  $w = [v,u]$
with $v \in I$ and $u \in \Gamma_{c}(TM)$.
\end{Lemma}
\proof
Choose $v \in I$ with $v(q) \neq 0$. There exist local co-ordinates
$x_1, \ldots, x_n$ 
and an open neighbourhood $W$ of $q$ such that $v|_{W} = \del_1$.
Choose $W$ to be a block centered around $q$, and
nest two smaller blocks (in local co-ordinates) 
inside, so that
$$
W \supset V \supset U \,.
$$   
We take 
$W = (-\varepsilon,\varepsilon)^{n}$,
$V = (- { \frac{2}{3}}\varepsilon,\frac{2}{3}\varepsilon)^{n}$, and
$U = (- { \frac{1}{3}}\varepsilon,\frac{1}{3}\varepsilon)^{n}$. 
Choose a smooth function $g$ on $M$ such that $g|_U = x_1$ 
and $g|_{M\backslash V}=0$. 
Require also that $\del_i g(x) = 0$ for $i \neq 1$ and
$(x_2, \ldots ,x_n) \in (- { \frac{1}{3}}\varepsilon,\frac{1}{3}\varepsilon)^{n-1}$.
Define $h := \del_1 g$, and set 
$\tilde{v} := [v,g\del_1] = h \del_1$.
Then $\tilde{v}|_U = \del_1$, $\tilde{v}|_{M\backslash V} = 0$, 
and most importantly
$\tilde{v} \in I$.

Now let $w \in \Gamma_{c}(TM)$ with $\supp(w) \subset U$.
We will find a global vector field $u \in \Gamma_{c}(TM)$ that realizes 
$[\tilde{v},u]=w$. 

For $i \neq 1$, the $i^{\mathrm{th}}$ component of the above reads
$h \del_1 u^i = w^i$. 
Set $u^i (x_1, \ldots , x_n) = 
\int_{-\infty}^{x_1} w^i(t,x_2 \ldots,x_n) \mathrm{d}t$ 
for $x_1 \leq \frac{1}{3}\varepsilon$, so that on $U$, where 
$\tilde{v} = \del_1$, we have $\del_1 u^i = w^i$
and thus $h \del_1 u^i = w^i$. This is obviously also correct for  
points outside $U$ with $x_1 \leq \frac{1}{3}\varepsilon$, as both $u^i$
and $w^i$ are zero.

For $\frac{1}{3}\varepsilon \leq x_1 \leq \frac{2}{3}\varepsilon$,
where $w=0$, 
we let $u^i (x_1, x_2, \ldots , x_n)$ be a constant function of
$x_1$, so that 
$u^i (x_1, x_2, \ldots , x_n) 
= u^i (\frac{1}{3} \varepsilon, x_2, \ldots , x_n)$.
We then have $\del_1 u^i (x_1, \ldots , x_n) = 0$, 
guaranteeing $h \del_1 u^i = w^i$.
Note that this does not effect the smoothness of $u^i$.

Finally, for $\frac{2}{3}\varepsilon \leq x_1 \leq \varepsilon$,
let $u^i$ tend to zero, and let $u^i$ be zero for 
$x_1 \geq \varepsilon$.
This can be done in such a way that $u^i$ remains smooth. 
Since both $h$ and $w$ are zero, we have $h \del_1 u^i = w^i$ on all of $M$.

The case $i = 1$ is handled similarly, the only difference being that 
the first component of 
$[\tilde{v},u]=w$ is now $h \del_1 u^1 - \sum_{j}u^j \del_j h = w^1$,
the term $u^1 \del_1 h$ of which cannot be dispensed with.
We have arranged for
$u^i(x)$, with $i \neq 0$, to be zero if   
$(x_2, \ldots ,x_n) \notin (- {
\frac{1}{3}}\varepsilon,\frac{1}{3}\varepsilon)^{n-1}$,
so that $\del_i h$ equals zero if $u^i$ does not.
This leaves us with 
$h \del_1 u^1 - u^1 \del_1 h = w^1$.

For $x_1 \leq \frac{1}{3}\varepsilon$,
one once again sets $u^1 (x_1, \ldots , x_n) = 
\int_{-\infty}^{x_1} w^1(t,x_2 \ldots,x_n) \mathrm{d}t$.
For $\frac{1}{3}\varepsilon \leq x_1 \leq \frac{2}{3}\varepsilon$ however,
one now has to define $u^i (x_1, \ldots , x_n)$ 
$=$ $h(x_1, \ldots , x_n)$ $u^i (\frac{1}{3}\varepsilon, x_2, \ldots , x_n)$
in order for $h \del_1 u^1 - u^1 \del_1 h = w^1$ to hold.
Since this renders $u^1$ zero on the boundary of $V$, one is then free 
to define $u^i$ to be zero on $M\backslash V$.

Thus, if $\supp(w) \subset U$, we see that 
$w = [\tilde{v},u]$ with $\tilde{v} \in I$ and
$u \in \Gamma_{c}(TM)$. This concludes the proof 
lemma \ref{techniek}, and thereby that of lemma \ref{Ranzig}. \qed

A maximal subalgebra $A$ of a Lie algebra $L$ is either selfnormalising 
or ideal. 
Indeed it is contained in its normaliser,
which therefore equals either $A$ or $L$. 

A theorem of Barnes \cite{Ba} states that a finite-dimensional 
Lie algebra is nilpotent if and only
if\footnote{
Actually, the `only if' part in Barnes' theorem is not written down in
\cite{Ba}, but this is immediately clear from the proof of Engel's theorem.
(See e.g. \cite{Hu}).
} every maximal subalgebra is an ideal. 
On the other extreme: 
\begin{Proposition}\label{perfectisideaal}
Let $L$ be a Lie algebra over a field $\mathbf{K}$, and
let 
$\cal{S}$ be the set of subspaces $A \subset L$ such that $A$ is
both an ideal and a maximal subalgebra.
Then
$$
[L,L] = \bigcap_{A \in \cal{S}} A \,.
$$ 
In particular, since the r.h.s reads `$L$' if $\cal{S} = \emptyset$,
$L$ is 
perfect $([L,L] = L)$
if and only if every maximal subalgebra is selfnormalising.
\end{Proposition}
\proof
Let $X \notin [L,L]$. Then choose 
$[L,L] \subseteq A \subsetneq L$ where $A$ has codimension 1 in $L$,
and $X \notin A$. A is an ideal maximal subalgebra,  
which does not contain $X$.
Thus $\bigcap_{A \in \cal{S}} A \subseteq [L,L]$.

Let $A$ be an ideal maximal subalgebra, and $X \notin A$.
Then $A + \mathbf{K} X$ is a subalgebra strictly containing $A$,
so that it must equal $L$. Thus $[L,L] = [A + \mathbf{K} X,A + \mathbf{K} X] 
\subseteq A$,
whence $[L,L] \subseteq \bigcap_{A \in \cal{S}} A$. \qed
 
As a corollary, we have the following well known statement: 
\begin{Corollary} \label{pstorm} The Lie algebra $\Gamma_{c}(TM)$ is perfect; 
$$[\Gamma_{c}(TM),\Gamma_{c}(TM)] = \Gamma_{c}(TM)\,.$$
\end{Corollary}
\proof 
According to lemma \ref{Ranzig}, the maximal ideals are precisely the ideals 
$I_{q}$ of vector fields in $\Gamma_{c}(TM)$ which are zero and
flat at $q$. 
$I_{q}$ is strictly contained in the subalgebra 
$A_{q}$ of vector fields which are zero at $q$, so that
no ideal is a maximal subalgebra. So every maximal subalgebra
is selfnormalising. \qed

\subsection{The splitting as a local map}\label{localmap}
With the main technical obstacles out of the way, we turn our attention to 
the sequence (\ref{rijtje}).
We now prove that $\sigma$ is a local map.

%


\begin{Lemma}\label{zwabbernoot}
Let $P \rightarrow M$ be a principal $G$-bundle over $M$, with $G$
any Lie group.
Let $\sigma: \Gamma_{c}(TM) \rightarrow \Gamma_{c}(TP)^{G}$
be a Lie algebra homomorphism splitting the exact sequence of Lie-algebras
\begin{equation}\label{srijtje}
0 \rightarrow \Gamma_{c}(TP)^{G}_{v} \rightarrow \Gamma_{c}(TP)^{G}
\rightarrow \Gamma_{c}(TM) \rightarrow 0 \,.
\end{equation}
Then $\sigma$ is local in the sense that 
$\pi (\supp(\sigma(v))) \subseteq \supp(v)$.
\end{Lemma}
\proof 
Any principal fibre bundle possesses an equivariant connection 1-form 
$\omega \in \Omega^{1}(P,\mathfrak{g})$, which
enables one to lift vector fields.
More precisely, the lifting map 
$\gamma :\Gamma_{c}(TM) \rightarrow \Gamma_{c}(TP)^{G}$ is defined 
by the requirement that 
$\omega_{p}(\gamma(v_{\pi(p)}))$ equal zero for all $p \in P$,
and that $\pi_{*} \circ \gamma$ be the identity. 
The lifting map $\gamma$ splits the exact sequence (\ref{srijtje}) as a sequence of  
$C^{\infty}(M)$-modules, but generally \emph{not} as a sequence 
of Lie algebras, since it need not be 
a homomorphism.


Define $f := \sigma - \gamma$. 
Then $\pi_{*} \circ f = 0$, 
so that $f$ is a map $\Gamma_{c}(TM) \rightarrow \Gamma_{c}(TP)^{G}_{v}$.
Since $\gamma$ is local by definition, 
the lemma will follow if we show
$f$ to be local.

The fact that $\sigma$ is a homomorphism, 
$[\sigma(v),\sigma(w)] - \sigma([v,w]) = 0$,
translates to 
\begin{equation}\label{vlijmenfileermes}
f([v,w]) - [f(v) , f(w)] =
[\gamma(v) , f(w)] - [\gamma(w), f(v)] 
+ 
[\gamma(v),\gamma(w)] - \gamma([v,w])
.
\end{equation}
The $G$-action effects a Lie algebra isomorphism 
$\Gamma_{c}(TP)^{G}_{v} \simeq \Gamma_{c}(\mathrm{ad}(P))$. 
From that perspective, we have for $s \in \Gamma_{c}(\mathrm{ad}(P))$ that
$[\gamma(v) , s] = \nabla_{v} s$, the covariant derivative
along $v$. The curvature of the connection is then given by
$R(v,w) = [\gamma(v),\gamma(w)] - \gamma([v,w])$.
We rewrite equation (\ref{vlijmenfileermes}) as
\begin{equation}\label{snelfornuis}
f([v,w]) - [f(v) , f(w)] = 
\nabla_v f(w) - \nabla_{w} f(v) + R(v,w)  \,.
\end{equation}
Pick $m \in M$, and identify the fibre of $\mathrm{ad}(P)$ 
over $m$ with $\mathfrak{g}$.
The restriction
$r_m : \Gamma_{c}(\mathrm{ad}(P)) \rightarrow \mathfrak{g}$
is a homomorphism of Lie algebras.
Now define $A_m := \{v \in \Gamma_{c}(TM) \,|\, v(m) = 0 \}$ to be the 
maximal subalgebra of vector fields that vanish at $m$,
and consider the map
$\tilde{f}_{m} : A_m \rightarrow \mathfrak{g}$
defined by $\tilde{f}_{m} = r_m \circ f |_{A_m}$.

Equation (\ref{snelfornuis}) transforms into 
\begin{equation}\label{cofflofot}
\tilde{f}_m([v,w]) - [\tilde{f}_m(v) , \tilde{f}_m(w)] = 
\nabla_v (\tilde{f}(w))_m - \nabla_{w} (\tilde{f}(v))_m + R(v,w)_m  \,,
\end{equation}
the r.h.s. of which vanishes
because $\nabla$ and $R$ are 
$C^{\infty}(M)$-linear in $v$ and $w$, 
and $v(m) = w(m) = 0$.
This means that $\tilde{f}_m$ is a Lie algebra homomorphism. 

Let us restrict $\tilde{f}_m$ even further to the (non-maximal) ideal 
$\Gamma_{c}(T(M-\{m\})) \subset A_m$, and note that
$\hat{f}_m : \Gamma_{c}(T(M-\{m\})) \rightarrow \mathfrak{g}$
is a homomorphism.
%
Its kernel $\mathrm{Ker}(\hat{f}_m)$ is therefore an ideal
in $\Gamma_{c}(T(M-\{m\}))$, and one of finite codimension at that. 
Indeed,
$$\Gamma_{c}(T(M-\{m\})) / \mathrm{Ker}(\hat{f}_m) 
\simeq \mathrm{Im}(\hat{f}_m) \subseteq \mathfrak{g}\,.
$$ 
According to lemma \ref{Ranzig} however, all proper ideals
are of infinite codimension, forcing 
$\mathrm{Ker}(\hat{f}_m) = \Gamma_{c}(T(M-\{m\}))$. 

But the vanishing of $\hat{f}_m$ for all $m \in M$ is tantamount to
locality of $f$ in the sense that $\supp(f(v)) \subseteq \supp(v)$. \qed

For the convenience of the reader, we gather some definitions made in the course 
of the proof. If $\sigma : \Gamma_{c}(TM) \rightarrow \Gamma_{c}(TP)^{G}$
is a Lie algebra homomorphism splitting $\pi_{*}$ and 
$\gamma : \Gamma_{c}(TM) \rightarrow \Gamma_{c}(TP)^{G}$ is the lift induced by
an equivariant connection, then we define $f := \sigma - \gamma$, which maps to 
the vertical vector fields, which we identify with $\Gamma(\mathrm{ad}(P))$.
We define $A_{m}$ to be the maximal subalgebra of vector fields which vanish at
$m$, and $\tilde{f}_{m} : A_{m} \rightarrow \mathfrak{g}$ to be the restriction
of $f$ to $A_{m}$, followed by the map 
$\Gamma(\mathrm{ad}(P)) \rightarrow \mathfrak{g}$ which picks out the fibre over 
$m$ and identifies it with $\mathfrak{g}$. 

\subsection{The splitting as a differential operator}\label{diffop}

In this section, we will prove that $\sigma$ is a differential operator
of finite order. 
Since $f$ is local, it defines a map from the sheaf of smooth sections 
of $TM$ to the sheaf of smooth sections of $\mathrm{ad}(P)$.
An elegant theorem of Peetre (\cite{Pe}) then says that $f$, and therefore 
$\sigma$, must be
a differential operator of locally finite order. 
All we need to do is find a global bound on the order, which will occupy
us for 
the remainder of the section.

The fact that $f$ is a differential operator of locally finite order 
means that for each $m \in M$, there exists an $r$ such that 
$H^{r}_m = \{v \in A_m \,|\, j^{r}_{m}(v) = 0\}$ is contained in 
$\ker(\tilde{f}_m)$, where $j^{r}_{m}(v)$ is the $r$-jet of $v$ at $m$.
Consequently, $\tilde{f}_{m}$ factors through the jet Lie algebra
$J^{r,0}_m(TM) := A_{m}/H^{r}_{m}$. Since we do not know $r$, we define 
$$\vecn := \lim_{\longleftarrow} A_m / H^{r}_{m}\,,$$
and remark that $\tilde{f}_m$ induces a homomorphism 
$\vecn \rightarrow \mathfrak{g}$. 
The Lie algebra $\vecn$ depends on $M$
only through its dimension $n$. 


Local co-ordinates provide one with a basis $x^{\vec{\alpha}}\del_i$,
where $x^{\vec{\alpha}}$ is shorthand for
$x_1^{\alpha_{1}} \ldots x_{n}^{\alpha_{n}}$.  
With
$\vecn^{k} = \mathrm{Span}\{ x^{\vec{\alpha}}\del_i \,|\, 
|\vec{\alpha}| = k +1 \, ,\, i = 1 \ldots n\} $, 
one may write
$$
\vecn = \bigoplus_{k=0}^{\infty}\vecn^{k}\,.
$$
That is, each element of $\vecn$ can be uniquely written as a 
\emph{finite} sum of homogeneous vector fields.
Note that $\vecn^{k}$ is the $k$-eigenspace of 
the Euler vector field $E := \sum_{i=1}^{n} x^i \del_i$. 
If $I$ is an ideal containing $v = \sum_{k=0}^{N} v_{k}$, 
then $\mathrm{ad}(E)^{j}v  = \sum_{k=0}^{N} k^{j} v_{k} \in I$ for all 
$j$, so that $v_{k} \in I$.
Thus any ideal splits into homogeneous components
$$
I  = \bigoplus_{k=0}^{\infty} I^{k}
$$ 
with $I^{k} = I \cap \vecn^{k}$. This renders the ideal structure 
of $\vecn$ more or less tractable, so that we may prove the following 
bound on the order of $\tilde{f}_m$.

\begin{Lemma}\label{hallekal}
The order of the differential operator $\sigma$ is at most 
$\mathrm{dim}(\mathfrak{g})$ unless $\mathrm{dim}(M) = 1$
and $\mathrm{dim}(\mathfrak{g}) = 2$, in which case the order is 
at most 3.
\end{Lemma}
\proof We closely follow Epstein and Thurston, \cite{ET}.
One checks by hand that the
only ideals of
$\mathrm{Vec}_{1} = \mathrm{Span}\{x^{k} \del \,|\, k \geq 1\}$
are 
$\mathrm{Span}\{x^2\del , x^k \del \,| \, k \geq 4 \}$
and $\mathrm{Span}\{x^k \del \,| \, k \geq N \}$ with $N \geq 1$.

Consider
$\mathrm{Vec}_{1}$ as a subalgebra of $\vecn$,
define $K$ to be the kernel of $\tilde{f}_m$, and let 
$K_1 := K \cap \mathrm{Vec}_{1}$. 
We then have injective homomorphisms
$$
\mathrm{Vec}_{1}/K_1 \hookrightarrow \vecn/K \hookrightarrow \mathfrak{g}\,,
$$
so that $\mathrm{dim}(\mathrm{Vec}_{1}/K_1) \leq \mathrm{dim}(\mathfrak{g})$.
As $K_1$ is an ideal, it must be of the shape mentioned above.
This leads us to conclude that $x_1^{k} \del_1 \in K$ 
for all $k > \mathrm{dim}(\mathfrak{g})$ unless 
$\mathrm{dim}(\mathfrak{g}) = 2$, in which case 
$x_1^{k} \del_1 \in K$ for all $k > 3$, and $x_1^2\del_1 \in K$.

The following short calculation shows that if $\mathrm{dim}(\mathfrak{g}) = 2$ and
$\mathrm{dim}(M) > 1$, then also $x_1^3 \del_1 \in K$.
As $K$ contains $x_1^2\del_1$, it also contains 
$[x_1^2\del_1 , x_1 \del_2] = x_1^2 \del_2$,
and thus $[x_1^2 \del_2 , x_1 x_2\del_1] = x_1^3 \del_1 - 2 x_1^2 x_2 \del_2$.
But by bracketing with $x_1^2 \del_2$ and $x_2 \del_1$ respectively,
we see that $x_1^3 \del_1 - 3 x_1^2 x_2 \del_2$ is in $K$,
ergo $x_1^3 \del_1 \in K$. 

The next step is to show that if $x_1^{s}\del_1 \in K$, then 
$K$ also contains all $x^{\vec{\alpha}} \del_i$ with $|\vec{\alpha}| = s$.
First of all, we remain in $K$ if we repeatedly apply 
$\mathrm{ad}(x_i \del_1)$ to $x_1^{s}\del_1$, to the effect of 
replacing $x_1$ by $x_i$ up to a nonzero factor.
This shows that $x^{\vec{\alpha}}\del_1 \in K$.
Then the relation 
$
x^{\vec{\alpha}}\del_i =
[x^{\vec{\alpha}}\del_1 , x_1 \del_i]  + x_1 \del_i x^{\vec{\alpha}} \del_1 
$
transfers membership of $K$ from right to left.

In the generic case $\mathrm{dim}(\mathfrak{g}) \neq 2$, 
$\mathrm{dim}(M) \neq 1$, we may conclude that  
$H_m^{\mathrm{dim}(\mathfrak{g})} \subset K$, so that the order of $\sigma$
is at most $\mathrm{dim}(\mathfrak{g})$.
In the exceptional case $\mathrm{dim}(\mathfrak{g}) = 2$,
$\mathrm{dim}(M) = 1$, the order of $\sigma$ is at most 3.\qed

In particular, $\sigma$ is a differential operator of finite rather
than locally finite order. 
Let us summarise our progress so far.
\begin{Proposition} \label{alhetvet}
Let $P$ be an infinitesimally natural principal $G$-bundle.
Then $\sigma : \Gamma_{c}(TM) \rightarrow \Gamma_{c}(TP)^{G}$ 
factors through the bundle of $k$-jets, where 
$k=3$ if $\mathrm{dim}(M) = 1$, $\dg = 2$ and
$k = \dg$ otherwise.  
If we identify 
equivariant vector fields on $P$ with
sections of 
the Atiyah bundle $TP/G$, we can therefore
define a bundle map 
$\nabla : J^{k}(TM)\rightarrow TP/G$ by
$\nabla (j_{m}^{k}(v)) := \sigma(v)_m$.
It makes the following diagram commute.
\begin{center}
\begin{tikzpicture}
\pgfsetzvec{\pgfpoint{0.385cm}{-0.385cm}}
\node (links) at (0, 0, 0) {$\Gamma(TM)$};
\node (rechts) at (2.5, 0, 0) {$\Gamma(TP/G)$};
\node (onder) at (1.3, -1.3, 0) {$\Gamma(J^{k}(TM))$};
\node (nab) at (2.2,-0.8,0) {$\nabla$};
\draw [->] (links) -- node[above] {$\sigma$}(rechts);
\draw [->] (links) -- (onder);
\draw [->] (onder) -- (rechts);
\end{tikzpicture}
\end{center}
\end{Proposition}
The point is that although $\sigma$ is defined only on sections, 
$\nabla$ comes from a veritable bundle map 
$J^{k}(TM) \rightarrow TP/G$.  
Note that although $\sigma$ was only defined on $\Gamma_{c}(TM)$, 
it extends to $\Gamma(TM)$ by locality. 

\section{Lie groupoids and algebroids of Jets}\label{drie}
The bundles $J^{k}(TM)$ and $TP/G$ are Lie algebroids, and it will be 
essential for us to prove that $\nabla : J^{k}(TM) \rightarrow TP/G$  
is a homomorphism of Lie algebroids. In order to do this, we will
first have a closer look at $J^{k}(TM)$ and $TP/G$, and at their 
corresponding Lie groupoids.

Let us first set some notation.
The jet group $G^{k}_{0,0}( \mathbf{R}^{n})$ is the group of 
$k$-jets of diffeomorphisms 
of $\mathbf{R}^n$ that fix $0$. It is the semi-direct 
product of $\gln$
and the connected, simply connected, unipotent Lie group of $k$-jets
that equal the identity to first order.

The subgroup $G^{+ k}_{0,0}(\mathbf{R}^{n})$ of orientation preserving 
$k$-jets is  
connected, but not simply connected. 
As $G^{+ k}_{0,0}(\mathbf{R}^{n})$ retracts to $\mathrm{SO}(\mathbf{R}^n)$, 
its homotopy group is isomorphic to 
$\{1\}$ if $n=1$, to
$\mathbf{Z}$ if $n=2$, and to
$\mathbf{Z} / 2\mathbf{Z}$ if
$n>2$. For brevity, we introduce the following notation.
\begin{Definition}
If $k > 0$, we denote $\pi_{1}(G^{+ k}_{0,0}(\mathbf{R}^{n}))$ by $Z$.
\end{Definition}
 
Thus for $n > 2$, the universal cover 
$\tilde{G}^{+ k}_{0,0}(\mathbf{R}^{n}) \rightarrow 
G^{+ k}_{0,0}(\mathbf{R}^{n})$ is $2 : 1$,
and restricts to the spin group over $\mathrm{SO}(\mathbf{R}^n)$.

%

\subsection{The Lie groupoid of $k$-jets}

In this section, we define the Lie groupoid $G^{k}(M)$ of $k$-jets, 
its maximal source-connected Lie subgroupoid $G^{+k}(M)$,
and the $k^{\mathrm{th}}$ order 
frame bundle $F^{k}(M)$.
   
Denote by $G^{k}_{m',m}(M)$ the manifold of $k$-jets at $m$ of diffeomorphisms
of $M$ which map $m$ to $m'$, and denote 
by $G^{k}(M) = \cup_{M \times M} G^{k}_{m',m}(M)$ the groupoid
of $k$-jets. If $j^{k}_{m}(\alpha)$ is a $k$-jet at $m$ of a 
diffeomorphism $\alpha$, 
then its source is 
$s(j^{k}_m(\alpha)) = m$, its target is $t(j^{k}_m(\alpha)) = \alpha(m)$,
and multiplication is given by concatenation.
  
Denote by $G^{k}_{*,m}(M)$ the manifold 
$s^{-1}(m)$ of $k$-jets with source $m$. The target map 
$t: G^{k}_{*,m}(M) \rightarrow M$ endows it with a structure of 
principal fibre bundle, the structure
group $G^{k}_{m,m}(M) \simeq G^{k}_{0,0}(\mathbf{R}^{n})$ acting 
freely and transitively on the right. 
As $G^{1}(M)_{*,m}$ is isomorphic to the frame bundle $F(M)$, 
one calls $G^{k}(M)_{*,m}$ the $k^{\mathrm{th}}$ order frame bundle, 
sometimes denoted $F^{k}(M)$.

\begin{Lemma}
Let $M$ be connected and let $k \geq 1$. Then $G^{k}(M)$ is
source-connected if and only if $M$ is not orientable.
\end{Lemma}
\proof
We may as well consider $k=1$, because the fibres of 
$G^{k}(M) \rightarrow G^{1}(M)$
can be contracted.
Each source fibre $G^{1}(M)_{*,m}$ of $G^{1}(M)$ is isomorphic to
the frame bundle. By definition, $M$ is oriented precisely when
the frames can be grouped 
into positively and negatively oriented ones.\qed

\begin{Definition} We define $G^{+k}(M)$ to be the maximal source-connected
Lie subgroupoid of $G^{k}(M)$, and denote its source fibre by 
$F^{+k}(M)$.
\end{Definition}

In the light of the previous lemma, this means that 
$G^{+ k}(M)$ is
the Lie groupoid of $k$-jets of
orientation preserving diffeomorphisms
if $M$ is orientable, and simply $G^{k}(M)$ if $M$ is not.

Note that the map $D : \mathrm{Diff}(M) \rightarrow \mathrm{Diff}(G^{k}(M))$
defined by
$D\alpha : j_m^{k}(\gamma) \mapsto j_m^{k}(\alpha \gamma)$
is a homomorphism of groups.
We will call it the $k^{\mathrm{th}}$ order derivative.
Because $D\alpha$ is
source-preserving 
and right invariant, 
it defines a homomorphism 
$\mathrm{Diff}(M) \rightarrow 
\mathrm{Aut}^{G^{k}_{m,m}(M)}(G^{k}_{*,m}(M))$, 
splitting the exact sequence of groups (\ref{flapperdrop}).
This makes $G_{*,m}^{k}(M) = F^{k}(M)$ into
a natural bundle.
Note that $F^{+k}(M)$ is infinitesimally natural.




\subsection{The Lie algebroid of $k$-jets}
The bundle $J^{k}(TM)$ possesses a structure of Lie algebroid, 
induced by the
Lie groupoid $G^{k}(M)$. 
We now describe
the Lie bracket on $\Gamma(J^{k}(TM))$ explicitly.
Later, in section \ref{loidhom}, we will use this to  
show that $\nabla$
is a Lie algebroid homomorphism. 

The Lie algebroid of $G^{k}(M)$ is a vector bundle
$A \rightarrow M$.
Its fibre 
$A_m$ is by definition the subspace of
the tangent space of $G^{k}(M)$ at $j_{m}^{k}(\id)$ 
which is annihilated by $ds$.
Sections of $A$ therefore correspond 
to right-invariant vector fields on $G^{k}(M)$ parallel to the source fibres.

Each curve in $G^{k}_{*,m}(M)$ through $j^{k}_{m}(\id)$ takes the shape
$c(t) = j^{k}_{m}(\alpha_{t})$ with $\alpha_{0} = \id$,
so that its tangent vector $a\in A_{m}$ takes the shape $a = j_{m}^{k}(v)$, 
with 
$v = {\del_{t}}|_{0} \alpha_{t}$. 
This shows that $A \simeq J^{k}(TM)$. 

The anchor $d t : J^{k}(TM) \rightarrow TM$ is easily seen to be the canonical
projection, so we shall denote it by $\pi$. The Lie bracket on 
$\Gamma(J^{k}(TM))$ however, 
which is defined as the restriction of 
the commutator bracket on $\Gamma(TG^{k}(M))$ to the right invariant 
source preserving vector fields, 
perhaps deserves 
some comment.


Define $J^{k,0}(TM)$ to be the kernel of $\pi$, and consider
the exact sequence of Lie algebras
$$
0 \rightarrow 
\Gamma(J^{k,0}(TM)) 
\rightarrow \Gamma(J^{k}(TM)) 
\stackrel{\pi}{\rightarrow} \Gamma(TM) 
\rightarrow 0\,.
$$ 
It is split by 
$j^{k} : 
\Gamma(TM) \rightarrow 
\Gamma(J^{k}(TM))
$, the infinitesimal version of the $k^{\mathrm{th}}$ order derivative. 
(The sequence of Lie algebroids of course
does not split, as differentiation is not linear over 
$C^{\infty}(M)$.)

We will describe the Lie bracket on 
$\Gamma(J^{k}(TM))$ by giving it on $\Gamma(TM)$
and $\Gamma(J^{k,0}(TM))$ separately, and then giving the action of 
$\Gamma(TM)$
on $\Gamma(J^{k,0}(TM))$.

\begin{Proposition}
Let $u$ and $u'$ be sections of $TM$, and 
let
$\tau : m \mapsto j_m^{k}(v_m)$ and 
$\tau' : m \mapsto j_m^{k}(v'_m)$ be sections of $J^{k,0}(TM)$. 
Then 
\begin{eqnarray*}
{}[j^{k}(u) , j^{k}(u')]_{m} &=& j_{m}^{k}([u,u'])\,,\\
{}[\tau , \tau']_{m} &=& j^{k}_{m}([v_{m},v_{m}'])\,,\\
{}[j^{k}(u) , \tau]_{m} &=& j^{k}_{m}([u, v_{m}]) + 
j_{m}^{k}( d_{u}|_{m}(x \mapsto v_x) )\,,
\end{eqnarray*}
where 
$d_{u}|_{m}(x \mapsto v_x)$ 
is the ordinary derivative at $m$ along $u$ of a map 
$M \rightarrow \Gamma(TM)$. 
Although both terms on the right hand side depend on the choice of 
$m \mapsto v_{m}$, their sum does not.
\end{Proposition}
\proof
The first equality is clear, as $j^{k}$ is a homomorphism of Lie algebras.
The second equality can be seen as follows. 
Consider the bundle of groups
$G^{k}(M)_{*,*} := \{j^{k}_{m}(\alpha) \in G^{k}(M)\, | \,\alpha(m) = m \}$,
with bundle map $s = t$. 
Its sections $\Gamma(G^{k}(M)_{*,*})$ form a group under 
pointwise multiplication,
the Lie algebra of which is $\Gamma(J^{k,0}(TM))$,
with the pointwise bracket.
As $\Gamma(G^{k}(M)_{*,*})$ acts from the left on 
$G^{k}(M)$ by 
$j^{k}_m(\gamma) \mapsto j_{m}^{k}(\alpha_m) \circ j_{m}^{k}(\gamma)$, 
respecting
both the source map and right multiplication, the inclusion
$\Gamma(J^{k,0}(TM)) \rightarrow \Gamma(J^{k}(TM))$ is a 
homomorphism of Lie algebras.
This proves the second line.

To verify the third line, we must choose a smooth
map $x \mapsto v_{x}$ from $M$ to $\Gamma(TM)$ such that
$\tau_{x} = j_{x}^{k}(v_{x})$ in a neighbourhood of $m$.
Each $v_{x}$ necessarily has a zero at $x$.  
If we denote $m(s) := \exp(su)m$, then the bracket 
$[j^{k}(u) , \tau]_{m}$ is by 
definition\footnote{The 
Lie algebra of the diffeomorphism group is
the Lie algebra of vector fields,
but the exponential map is given by
$v \mapsto \exp(-v)$, where $\exp$ denotes the
unit flow along $v$.
This is why the groupoid commutator might seem odd at first sight.} 
minus the mixed second derivative along $s$ and $t$ at $0$ of
the groupoid commutator
$$
j^{k}_{m(s)}(\exp(-su)) 
j^{k}_{m(s)}(\exp(-tv_{m(s)})) 
j^{k}_{m}(\exp(su)) 
j_{m}^{k}(\exp(t v_{m}))
\,,$$
which is just 
$
j^{k}_{m} \left(
\exp(-su) 
\exp(-tv_{m(s)}) 
\exp(su) 
\exp(t v_{m})
\right)\,.
$
The terms not involving derivatives of $s \mapsto m(s)$ yield
$j_{m}^{k}([u,v_{m}])$, and the terms which do provide the extra
$j^{k}_{m}(d_{u}|_{m}(x \mapsto v_{x}))$. \qed

\subsection{The Gauge groupoid and its algebroid}
Given a principal $G$-bundle $\pi: P \rightarrow M$, one can define the
gauge groupoid $(P \times P)/G$, that is the pair groupoid modded out by the
diagonal action. 
Source and target 
come from projection on the second and first term respectively, 
$M \hookrightarrow (P \times P)/G$ as $\id_{\pi(p)} = [(p,p)]$, and 
multiplication is well defined by $[(r,q)] \circ [(q,p)] = [(r,p)]$.
An element $[(q,p)]$ corresponds precisely to a $G$-equivariant 
diffeomorphism $\pi^{-1}(p) \rightarrow \pi^{-1}(q)$, and the product to
concatenation of maps.

Its Lie algebroid $TP/G$ is sometimes called the Atiyah algebroid.
Indeed, 
the subspace of $T_{\id_{m}}((P \times P)/G)$ which annihilates $ds$
is canonically $(TP/G)_{m}$.
A section of $TP/G$ can be identified with a $G$-equivariant
section of $TP$, endowing $\Gamma(TP/G)$ with the Lie bracket
that comes from $\Gamma(TP)^{G}$.

\subsection{The splitting as a homomorphism of Lie algebroids}\label{loidhom}

The point of considering the Lie algebroid structure of $J^{k}(TM)$
was of course to prove the following.
\begin{Lemma}
The map $\nabla : J^{k}(TM) \rightarrow TP/G$ is a homomorphism of
Lie algebroids.
\end{Lemma}
 
\proof The fact that $\nabla$ respects the anchor is immediate.
As $\pi_{*} \circ \nabla (j_{m}^{k}(v_m)) = 
\pi_{*} \circ \sigma(v_{m})(m) = v_{m}(m)$, it equals 
$\pi \circ j^{k}(v_{m})$ in the point $m$.

We now show that $\nabla : \Gamma(J^{k}(TM)) \rightarrow \Gamma(TP/G)$
is a homomorphism of Lie algebras. First of all, the restriction 
of $\nabla$ to $j^{k}(\Gamma(TM))$ is a homomorphism. Indeed,
$[\nabla(j^{k}(v)) , \nabla(j^{k}(w))] = [\sigma(v), \sigma(w)]$,
which equals $\sigma([v,w])$ because $\sigma$ is a homomorphism.
This in turn is just $\nabla(j^{k}([v,w]))$, so that
$[\nabla(j^{k}(v)) , \nabla(j^{k}(w))] = \nabla([j^{k}(v),j^{k}(w)] )$.

Secondly, its restriction to $\Gamma(J^{k,0}(TM))$ is a homomorphism.
If $\tau_x = j^{k}_{x}(v_x)$ and $\upsilon_x = j^{k}_{x}(w_x)$ 
are sections of $J^{k,0}(TM)$, then 
$\nabla \tau$ and $\nabla \upsilon$ are in the kernel of the anchor.
This implies that
their commutator at a certain point $m$ depends only on their
values at $m$, not on their derivatives.
To find the commutator at $m$, we may therefore
replace $j^{k}_{x}(v_x)$ by $j^{k}_{x}(v_{m})$ and likewise 
$j^{k}_{x}(w_x)$ by $j^{k}_{x}(w_{m})$. We then see that
$[\nabla j^{k}_{x}(v_x) , \nabla j^{k}_{x}(w_x)]_{m}$
$=$ $[\nabla j^{k}_{x}(v_{m}) , \nabla j^{k}_{x}(w_{m})]_{m}$.
We already know that this is
$\nabla( j^{k}_{x}([v_{m},w_{m}]))_{m}$,
so that
$[\nabla(\tau),\nabla(\upsilon)]_{m} = \nabla([\tau,\upsilon])_{m}$. 

The last step is to show that $\nabla$ respects the bracket between 
$j^{k}(\Gamma(TM))$ and $\Gamma(J^{k,0}(TM))$.
Again, let
$j^{k}(v)$ be an element of the former and $\tau_x = j^{k}_{x}(w_x)$ 
of the latter.
Considered as an equivariant vector field on $P$, the vertical 
vector field $\nabla(\tau)$ takes the value $\sigma_{p}(w_{\pi(p)})$
at $p \in P$.
Then $[\nabla(j^{k}(v)) , \nabla(\tau)]$ is the Lie derivative along 
$\sigma(v)$ of the vertical vector field $\sigma_{p}(w_{\pi(p)})$.
Differentiating along $\sigma(v)$ is done by considering
$(p,p') \mapsto \sigma_{p}(w_{\pi(p')})$, differentiating w.r.t. 
$p$ and $p'$ separately, and then putting $p=p'$.
This results in
$[\sigma(v), \sigma_{p}(w_{\pi(p)})]_{p_0} = 
[\sigma(v), \sigma(w_{\pi(p_0)})]_{p_0} + 
\sigma(d_{v}|_{\pi(p_0)} (x \mapsto w_{x}))$, 
which is in turn the same as 
$\sigma_{p_0}([v,w_{\pi(p_0)}] + d_{v}|_{\pi(p_0)} (m \mapsto w_{m})))$,
so that $[\nabla(j^{k}(v)) , \nabla(\tau)] = \nabla([j^{k}(v),\tau])$
as required. 
Therefore $\nabla$ must be a homomorphism on all of $\Gamma(J^{k}(TM))$. \qed

\begin{Definition}
A connection $\nabla$ of a Lie algebroid $A$ on a vector bundle $E$
is by definition a bundle map of $A$ into $\mathrm{DO}^{1}(E)$,
the first order differential operators on $E$, which respects the anchor.
If moreover it is a morphism of Lie algebroids, then the connection is 
called flat. A flat connection of $A$ on $E$ is also called a representation
of $A$ on $E$.
\end{Definition}

This explains our notation for the map $\nabla$ induced by $\sigma$.
Given a representation $V$ of $G$, one may form the 
associated vector bundle $E := P \times_{G} V$. The map $\nabla$
then defines a Lie algebroid homomorphism of $\Gamma(J^{k}(TM))$ 
into the Lie algebroid of first order differential operators on $E$.
(Simply consider a section of $E$ as a $G$-equivariant function
$P \rightarrow V$, and let $\Gamma^{G}(TP)$ act by Lie derivative.)
By definition, this is a flat connection, or equivalently a 
Lie algebroid representation.

\section{The classification theorem}\label{vier}
We use the fact that $\nabla : \Gamma(J^{k}(TM)) \rightarrow \Gamma(TP/G)$ 
is a homomorphism of Lie algebroids
to find a corresponding homomorphism of Lie groupoids.
This will give us the desired classification of 
infinitesimally natural principal fibre
bundles.

\subsection{Integrating a homomorphism of Lie algebroids}
The following theorem states that homomorphisms of Lie algebroids
induce homomorphisms of Lie groupoids if the initial groupoid is 
source-simply connected.
\begin{Theorem}[Lie II for algebroids]
Let $G$ and $H$ be Lie groupoids, with corresponding Lie algebroids 
$A$ and $B$ respectively. Let $\nabla : A \rightarrow B$
be a homomorphism of Lie algebroids. If $G$ is source-simply 
connected, then there exists a unique homomorphism
$G \rightarrow H$ of Lie groupoids which integrates $\nabla$.
\end{Theorem}
The result was probably announced first in \cite{Pr}, proofs have 
appeared e.g. in \cite{MX} and \cite{MM}.
We follow \cite{CF}, which the reader may consult for details.

{\bf Sketch of proof \,}
The idea is that $\nabla$ allows one to lift 
piecewise smooth paths 
of constant source in $G$ to piecewise smooth paths of constant source 
in $H$. Source-preserving piecewise smooth homotopies in $G$ of course 
do not affect the endpoint of the path in $H$, so that,
if $G$ is source-simply connected, 
one obtains a map $G \rightarrow H$ by identifying 
elements $g$ of $G$ with equivalence classes of source preserving
paths from $\id_{s(g)}$ to $g$. 
One checks that this is 
the unique homomorphism of Lie groupoids integrating $\nabla$. 
\qed

Unfortunately, $G^{k}(M)$ is not always source connected, 
let alone source-simply connected.
Recall that $G^{+k}(M)$ is the maximal source-connected Lie subgroupoid of 
$G^{k}(M)$, and therefore has the same Lie algebroid $J^{k}(TM)$. 

We define $\tilde{G}^{+k}(M)$ to be the 
set of piecewise smooth, source preserving paths in $G^{+k}(M)$
beginning at an identity, 
modulo piecewise smooth, source preserving homotopies.
It is a smooth manifold because $G^{+k}(M)$ is, and a Lie groupoid
under the unique structure making the projection on the endpoint
$\tilde{G}^{+k}(M) \rightarrow G^{+k}(M)$ into a morphism of 
groupoids.
Explicitly, the multiplication is given as follows. 
If $g(t)$ a path from $\id_{m}$ to $g(1)_{m'm}$,
and $h(t)$ a path from $\id_{m'}$ to $h(1)_{m''m'}$,
then the product $[h] \circ [g]$ is 
$[(h \cdot g(1)) * g]$, where the dot denotes groupoid multiplication 
and the star concatenation of paths. The proof of associativity is the usual
one. 

Note that the source fibre $\tilde{G}^{+k}(M)_{*,m}$ is precisely
the universal cover of the connected component of
the $k^{\mathrm{th}}$ order frame bundle 
$G^{+k}(M)_{*,m} = F^{+k}(M)$.
In order to cut down on the subscripts, we introduce new notation 
for $\tilde{G}^{+k}(M)_{*,m}$ and its structure group 
$\tilde{G}^{+k}(M)_{m,m}$.
\begin{Definition} 
We denote the universal cover of the connected component of the 
$k^{\mathrm{th}}$ order frame bundle by $\tf$, and
its structure group by $G(k,M)$.  
\end{Definition}
It is an infinitesimally natural bundle because 
$F^{+k}(M)$ is. 
Note that $G(k,M)$ is not the universal cover of 
$G_{m,m}^{+k}(M)$,
but rather its extension by $\pi_{1}(F^{k}(M))$. 
As $\pi_{1}(F^{k}(M)) = \pi_{1}(F(M))$, we have the exact 
sequence of groups
$$
1 \rightarrow \pi_{1}(F(M)) \rightarrow G(k,M) \rightarrow 
G_{m,m}^{+k}(M) \rightarrow 1\,.
$$
The group $G_{m,m}^{+k}(M)$ in turn is isomorphic to 
$G_{0,0}^{+k}(\mathbf{R}^{n})$ if $M$ is orientable, and to
$G_{0,0}^{k}(\mathbf{R}^{n})$ if it is not.

\subsection{Classification}\label{Classification}

Now that we've found a source-simply connected Lie groupoid with 
$J^{k}(TM)$ as Lie algebroid, we can finally apply 
Lie's second theorem for algebroids to obtain the following.
\begin{Proposition}\label{knolletje}
If $\sigma$ splits the exact sequence of Lie algebras (\ref{rijtje}),
then it induces a morphism of groupoids 
$\exp\nabla : \tilde{G}^{+k} \rightarrow (P \times P)/G$ such that the following diagram
commutes, with $\exp_{m}$ the flow along a vector field starting at 
$\id_{m}$.
\end{Proposition}
\begin{center}
\begin{tikzpicture}
\pgfsetzvec{\pgfpoint{0.385cm}{-0.385cm}}
\node (links) at (0, 0) {$\Gamma(TM)$};
\node (midden) at (3, 0) {$\Gamma(TP/G)$};
\node (onderlinks) at (1.5, -1.5) {$\Gamma(J^{k}(TM))$};
\node (rechts) at (6,0) {$(P \times P )/G$};
\node (onderrechts) at (4.5, -1.5) {$\tilde{G}^{+k}(M)$};
\node at (0.3,-0.8) {$j^{k}$};
\node at (4.7,-0.7) {$\exp\nabla$};
\node at (1.95,-0.7) {$\nabla$};
\draw [->] (links) -- node[above] {$\sigma$}(midden);
\draw [->] (links) -- (onderlinks);
\draw [->] (onderlinks) -- (midden);
\draw [->] (midden) -- node[above] {$\exp_{m}$} (rechts);
\draw [->] (onderlinks) -- node[above] {$\exp_{m}$} (onderrechts);
\draw [->] (onderrechts) -- (rechts);
\end{tikzpicture}
\end{center}
\proof As $\tilde{G}^{+k}(M)$ is a source-simply connected Lie groupoid with 
$J^{k}(TM)$ as Lie algebroid, we can apply 
Lie's second theorem for algebroids. \qed

It is perhaps worth wile to formulate this
for general transitive Lie groupoids.
\begin{Proposition}
Let $\mathcal{G}\rightrightarrows M$ be a transitive Lie groupoid, 
with Lie algebroid $A$. The kernel of the anchor $K$ is then a 
bundle of Lie algebras with fixed dimension $d$.  Suppose that the sequence 
$$
0 \rightarrow \Gamma(K) \rightarrow \Gamma(A) \rightarrow 
\Gamma(TM) \rightarrow 0\,,
$$
with $K$ the kernel of the anchor, splits as a sequence of Lie algebras. 
Then this splitting is induced by a morphism of Lie algebroids 
$\nabla : J^{k}(TM) \rightarrow A$,  
and there is a corresponding morphism of Lie groupoids 
$\tilde{G}^{+k}(M) \rightarrow \mathcal{G}$. The number $k$ is at most $3$
if $d$ is $2$ and $\mathrm{dim}(M)=1$, and at most $d$ otherwise.
\end{Proposition}
\proof (Or rather a flimsy sketch thereof.) 
Analogous to the case of the gauge groupoid. \qed

We have paved the way for a classification of
infinitesimally natural principal fibre bundles.
\begin{Theorem}\label{grondsmakreefknol}
Let $\pi : P \rightarrow M$ be an infinitesimally natural principal 
$G$-bundle with splitting $\sigma$ of (\ref{rijtje}).
Then there exists a group homomorphism 
$\rho : G(k,M) \rightarrow G$ such that 
the bundle $P$ is 
associated to $\tf$ through $\rho$, i.e. 
$$P \simeq \tf \times_{\rho} G\,.$$
Moreover, $\sigma$ is induced by the canonical one for
$\tf$.
\end{Theorem}
\proof Fix a base point $m$ on $M$. The map $\exp \nabla$ 
yields a homomorphism of groups
$\rho : \tilde{G}_{m,m}^{+k}(M) \rightarrow ((P\times P)/G)_{m,m}$, 
the latter isomorphic to $G$, 
the former to 
$G(k,M)$.

The map
$\tilde{G}_{*,m}^{+k}(M) \times_{\rho}  
((P\times P)/G)_{m,m} \rightarrow ((P\times P)/G)_{*,m}$
which is given by $(g_{m' , m} , p_{m,m}) \mapsto 
(\exp \nabla (g_{m' , m})) \cdot p_{m,m}$
is well defined and injective because two pairs share the same 
image if and only
if they are equivalent modulo $\tilde{G}_{m,m}^{+k}(M)$.
It is also surjective and $G$-equivariant, and hence an
isomorphism of principal $G$-bundles. As 
$((P \times P)/G)_{*,m} \simeq P$
and
$\tilde{G}_{*,m}^{+k}(M) \times_{\rho}  
((P\times P)/G)_{m,m} \simeq \tf \times_{\rho} G$,
the equivalence is proven. The remark on $\sigma$ 
follows from the construction.
\qed

This classifies the infinitesimally natural principal fibre bundles. 
They are all associated (via a group homomorphism) 
to the bundle $\tilde{G}^{+k}_{*,m} = \tf$.

The classification of natural principal fibre bundles is now an easy
corollary. The following well known result (\cite{PT}, \cite{Te}) states that they are precisely the
ones associated to $G^{k}_{*,m}(M) = F^{k}(M)$.
\begin{Corollary}\label{PalaisTerng}
Let $\pi : P \rightarrow M$ be a natural principal $G$-bundle
with local splitting $\Sigma$ of (\ref{flapperdrop}). 
Then $P$ is associated to 
$F^{k}(M)$.
That is, there exists a 
homomorphism $\rho : G^{k}_{0,0}(\mathbf{R}^{n}) \rightarrow G$ such that 
$$
P \simeq F^{k}(M) \times_{\rho} G\,.
$$
Moreover, $\Sigma$ is induced by the canonical one for $F^{k}(M)$.
\end{Corollary}
\proof
As the homomorphism 
$\Sigma : \mathrm{Diff}_{c}(M) \rightarrow \mathrm{Aut}_{c}(P)$
is local, it induces
a homomorphism of groupoids $\Sigma : 
\mathrm{Germ}(M) \rightarrow (P \times P)/G$,
with $\mathrm{Germ}(M)$ the groupoid of germs of diffeomorphisms of $M$.
We need but show that $\Sigma$ factors through 
$j^{k} : \mathrm{Germ}(M) \rightarrow G^{k}(M)$ for some $k > 0$,
cf. the proof of theorem \ref{grondsmakreefknol}.
 
The Lie algebra homomorphism 
$\sigma : \Gamma_{c}(TM) \rightarrow TP/G$ defined by
$\sigma(v) := \del_{t}|_{0} \Sigma(\exp(tv))$ 
is local by assumption, and according to proposition \ref{alhetvet} 
it factors through 
the $k$-jets for some $k>0$.
It suffices to show that $\Sigma(\phi)_{m,m} = \id_{m,m}$
for any $\phi \in \mathrm{Germ}_{m,m}(M)$ 
that agrees with the identity to $k^{\mathrm{th}}$ order at $m$.

In local co-ordinates $\{x^{i}\}$,
we write $\phi^{i}(x) = x^{i} + v^{i}(x)$, where 
$v : \mathbf{R}^{n} \rightarrow \mathbf{R}^{n}$
vanishes to $k^{\mathrm{th}}$ order.
We define the one parameter family of germs of diffeomorphisms 
$\phi^{i}_{t}(x) := x^{i} + t v^{i}(x)$.
Then 
$\del_{t}|_{\tau} \Sigma(\phi_{\tau})^{-1} \Sigma(\phi_{t})_{m} = 0$,
as it equals $\sigma_{m}( \del_{t}|_{\tau} \phi^{-1}_{\tau} \phi_{t})$,
the image of a vector field that 
vanishes to 
order $k$ at $m$. Therefore  
$t \mapsto \Sigma(\phi_{t})_{m,m}$ is constant, and
$\Sigma(\phi)_{m,m} = \id_{m,m}$ as required.
\qed 

To summarise: natural principal fibre bundles are associated to a
higher frame bundle, whereas infinitesimally natural 
principal fibre bundles are associated to the universal cover of 
a higher frame bundle.  



\subsection{The bundle $\tf$}
The above considerations prompt a few remarks on the 
universal cover of the connected component of the frame bundle
$\tf$, and on its (disconnected) structure group $G(k,M)$.
Recall that they are just the source fibre $\tilde{G}^{+k}_{*,m}(M)$
and isotropy group $\tilde{G}^{+k}_{m,m}(M)$ of $\tilde{G}^{+k}(M)$.
 
\subsubsection{General manifolds}
If $\pi_{1}(M)$ is the homotopy groupoid of $M$, define
the homomorphism of groupoids $\mathrm{Pr} : \pi_{1}(M)_{m',m} \rightarrow 
\pi_{0}(G_{m',m}(M))$
by lifting a path in $M$ to a path in $G^{k}(M)$ with fixed source,
and taking the connected component of its end point.
It makes
\begin{center}
\begin{tikzpicture}
\pgfsetzvec{\pgfpoint{0.385cm}{-0.385cm}}
\node (links) at (0, -0.7, 0) {$\tilde{G}_{m',m}^{+k}(M)$};
\node (rechts) at (2.8, 0, 0) {$\pi_{1}(M)_{m',m}$};
\node (linksonder) at (2.8, -1.4, 0) {$G_{m',m}^{k}(M)$};
\node (rechtsonder) at (5.8,-0.7, 0) {$\pi_{0}(G_{m',m}^{k}(M))$};
\draw [->] (links) -- (rechts);
\draw [->] (links) -- (linksonder);
\draw [->] (linksonder) -- (rechtsonder);
\draw [->] (rechts) -- (rechtsonder);
\end{tikzpicture}
\end{center}
into a commutative diagram.

Define 
$(G^{k}(M) \times \pi_{1}(M))^{\mathrm{Pr}}$
to be the groupoid of pairs $(g,[f])$
such that $\pi_{0}(g) = \mathrm{Pr}([f])$. 
If $M$ is orientable, this is simply 
$G^{+k}(M) \times \pi_{1}(M)$.
The map of groupoids
$\tilde{G}^{+k}(M) \rightarrow 
(G^{k}(M) \times \pi_{1}(M))^{\mathrm{Pr}}$
is well defined and surjective. 
It restricts to a covering map of principal fibre bundles
\begin{equation}\label{tauprojectie}
\tau : \tilde{G}_{*,m}^{+k}(M) \rightarrow 
(G^{k}(M) \times \pi_{1}(M))_{*,m}^{\mathrm{Pr}}\,.
\end{equation}
The kernel of the corresponding cover of groups 
is precisely
$i_{*}\pi_{1}(G_{m,m}^{+k}(M))$, with 
$i : G_{m,m}^{+k}(M) \rightarrow G_{*,m}^{+k}(M)$
the inclusion. 
Note that $i_{*}$ has a 
nonzero kernel precisely when a vertical loop is contractible
in $G_{*,m}^{+k}(M)$, but not by a homotopy
which stays inside the fibre. Denoting 
$\pi_{1}(G_{m,m}^{+k}(M))$ by $Z$, we obtain the exact sequence
\begin{equation}\label{centraalgroep}
1 \rightarrow Z/\mathrm{Ker}(i_{*}) 
\stackrel{i_{*}}{\rightarrow} 
\tilde{G}_{m,m}^{+k}(M)
\stackrel{\tau}{\rightarrow}
(G^{k}(M) \times \pi_{1}(M))_{m,m}^{\mathrm{Pr}}
\rightarrow
1\,.
\end{equation}
A moment's thought reveals that this extension is central: 
if $g(t)$ is a path in
$G_{m,m}^{+k}(M)$ and $h(t)$ one in 
$G_{*,m}^{+k}(M)$, then both $h * (i \circ g)$ and 
$(i \circ g)\cdot h(1) * h$
can be homotoped into $t \mapsto h(t)g(t)$.

We may as well restrict attention to the case $k=1$, 
in which $G_{*,m}^{1}(M)$ is the frame bundle $F(M)$.
Indeed, as 
$G^{+k}_{*,m}(M) \rightarrow G^{1,+}_{*,m}(M)$ has contractible fibres,
$\tilde{G}^{+k}_{*,m}(M)$ is just
the pullback of $G^{+k}_{*,m}(M)$ along 
$\tilde{G}_{*,m}^{1,+}(M) \rightarrow G_{*,m}^{1,+}(M)$.

\subsubsection{Orientable manifolds}
For orientable manifolds, the situation simplifies.
If we identify the connected component of $G_{m,m}^{1}(M)$ with $\glnp$, we 
obtain a homomorphism 
$i_{*}$ of $\tglnp$ into $\tilde{G}^{1,+}_{m,m}(M)$.
There is a second homomorphism
$\pi_{1}(F(M)) \rightarrow \tilde{G}^{1,+}_{m,m}(M)$.
Their images intersect in $Z/\mathrm{Ker}(i_{*})$, and commute
by an argument similar to the one on centrality of (\ref{centraalgroep}).
If we define $(\tglnp \times \pi_{1}(F(M)))_{Z}$ to be the 
quotient of $\tglnp \times \pi_{1}(F(M))$ by the equivalence
$(gz,h) \sim (g,zh)$,
we can regard it as a subgroup of $\tilde{G}^{1,+}_{m,m}(M)$.
Note that if $\mathrm{Ker}(i_{*})$ is nonzero, the above equivalence
relation sets it to 1.

If $M$ is orientable, we may restrict our attention to $F^{+}(M)$,
which has connected fibres. Any path in
$F^{+}(M)$ which starts and ends in the same fibre can therefore be obtained 
by combining a closed loop with a path in $\glnp$.
For orientable manifolds, we thus have 
$
\tilde{G}_{m,m}^{1,+}(M) \simeq (\tglnp \times \pi_{1}(F^{+}(M)))_{Z}
$, and in the same vein
\begin{equation}\label{orient}
G(k,M) \simeq 
(G(k,\mathbf{R}^{n}) \times \pi_{1}(F^{+}(M)))_{Z}\,.
\end{equation}

\subsubsection{Spin manifolds}\label{Spinmanifolds}
Let $M$ be an orientable manifold, equipped with a 
pseudo-Riemannian metric $g$ of signature 
$\eta \in \mathrm{Bil}(\mathbf{R}^{n})$.
Then $OF^{+}_{g} := \{f \in F^{+}(M) \,|\, f^* g = \eta\}$
is the bundle of positively oriented orthogonal frames. 
A spin structure is then by definition an
$\widetilde{\mathrm{SO}}(\eta)$-bundle\footnote{
There is a subtlety here. Suppose 
$\eta$ has indefinite signature, say $(3,1)$.
The group $\mathrm{SO}(3,1)$
has 2 connected components,  
so that a universal cover does not exist. 
As it is a subgroup of the simply connected group 
$\glnvp$,
we simply define
$\widetilde{\mathrm{SO}}(3,1)$ to be $\kappa^{-1}(\mathrm{SO}(3,1))$ 
with $\kappa : \tglnvp \rightarrow \tglnvp$ the covering map.
Thus $\widetilde{\mathrm{SO}}(3,1)$ is, perhaps surprisingly, 
not isomorphic to 
the 2-component spin group $\mathrm{Spin}(3,1)$. 
Indeed,
if
$T$ is time inversion and $P$
is the inversion of $3$ space co-ordinates, then
$(PT)^{2} = \one$
in $\widetilde{\mathrm{SO}}(3,1)$,
as opposed to
$(PT)^{2} = -\one$ in $\mathrm{Spin}(3,1)$.
Therefore 
$\pi^{-1}(\pm \one) \simeq 
\mathbf{Z}/2\mathbf{Z} \times \mathbf{Z}/2\mathbf{Z}$
in $\widetilde{\mathrm{SO}}(3,1)$,
whereas
$\pi^{-1}(\pm \one) = \mathbf{Z}/4\mathbf{Z}$ in $\mathrm{Spin}(3,1)$
(see \cite{BDWGK}).  

Of course 
the connected component of unity of $\widetilde{\mathrm{SO}}(3,1)$
and that of $\mathrm{Spin}(3,1)$ are both isomorphic 
to $\mathrm{SL}(\mathbf{C}^{2})$, so that
none of this is relevant
if $M$ is both orientable and time-orientable,
i.e. if the structure group of 
the frame bundle reduces to $\mathrm{SO}^{\uparrow}(3,1)$.
}
$Q$ over $M$, plus a map $u : Q \rightarrow OF^{+}_{g}$ 
such that
\begin{center}
\begin{tikzpicture}
\pgfsetzvec{\pgfpoint{0.385cm}{-0.385cm}}
\node (lg) at (-1,0.7) {$\widetilde{\mathrm{SO}}(\eta)$};
\node (rg) at (1,0.7) {$\mathrm{SO}(\eta)$};
\node  (rpijl) at (1,0.3) {$\curvearrowleft$};
\node  (lpijl) at (-1,0.3) {$\curvearrowleft$};
\node (lb) at (-1,0) {$Q$};
\node (rb) at (1,0) {$OF^{+}_{g}$};
\node (lo) at (0,-1.2) {$M$};
\draw [->] (lg) --node[above]{$\kappa$} (rg);
\draw [->] (lb) --node[above]{$u$} (rb);
\draw [->] (lb) -- (lo);
\draw [->] (rb) -- (lo);
\end{tikzpicture}
\end{center}
commutes, with $\kappa$ the canonical homomorphism 
$\widetilde{\mathrm{SO}}(\eta) \rightarrow \mathrm{SO}(\eta)$.
A manifold is called spin if it admits a spin structure.

Define $\hat{Q} := Q \times_{\widetilde{\mathrm{SO}}(\eta)} \tglnp$, 
and again denote the induced
map $\hat{Q} \rightarrow F^{+}(M)$ by $u$. 
As any cover of $F^{+}(M)$ by a $\tglnp$-bundle
can be obtained in this way, there is a 1:1 correspondence
between spin covers of $OF^{+}_{g}(M)$ and $F^{+}(M)$.
In particular, whether or not $M$ is spin does not depend on the metric.

The Serre spectral sequence gives rise to the exact sequence
\begin{equation} \label{piraterij}
1 \rightarrow Z/\mathrm{Ker}(i_{*}) \rightarrow 
\pi_{1}(F^{+}(M)) \rightarrow
\pi_{1}(M) \rightarrow 1\,.
\end{equation}
The following proposition is well known.
\begin{Proposition}\label{klopboor?}
A spin structure exists
if and only if $i_{*} : Z \rightarrow \pi_{1}(F^{+}(M))$
is injective and
(\ref{piraterij}) splits 
as a sequence of groups. If spin structures exist, 
then 
equivalence classes of spin covers correspond to
splittings of (\ref{piraterij}).
\end{Proposition}
\proof This will follow from theorem \ref{naturalspin} later
on, but
see e.g. \cite{Mo} for an independent proof. Our criterion for 
$M$ to be spin is equivalent to the vanishing of the second Stiefel-Whitney
class, see e.g. \cite{LM}.
\qed

{\bf Remark\quad}
In terms of group cohomology, one can consider the sequence 
(\ref{piraterij}) as an element
$[\omega] \in H^{2}(\pi_{1}(M) , Z/\mathrm{Ker}(i_*))$.
Spin bundles exist if and only if both $\mathrm{Ker}(i_*)$  
and $[\omega]$ are trivial, in which case they are indexed by
$H^{1}(\pi_{1}(M) , Z)$.

If a spin structure exists, then $\widetilde{F}^{+}$
is simply the pullback along
the universal cover $\tilde{M} \rightarrow M$ of  
$\hat{Q} \rightarrow M$. 
The picture then becomes 
\begin{center}
\begin{tikzpicture}
\pgfsetzvec{\pgfpoint{0.385cm}{-0.385cm}}
\node (mb) at (0,0) {$u_{*}\big( F^{+k}(M) \big)$};
\node (mo) at (0,-1.5) {$\hat{Q}$};
\node (lb) at (-3,0) {$\tilde{F}^{+k}(M)$};
\node (rb) at (3,0) {$F^{+k}(M)$};
\node (lo) at (-3,-1.5) {$\tilde{F}^{+}(M)$};
\node (ro) at (3,-1.5) {$F^{+}(M)$};
\node (moo) at (0,-3) {$M$};
\node (loo) at (-3,-3) {$\tilde{M}$};
\draw [->] (lb) -- (mb);
\draw [->] (mb) -- (rb);
\draw [->] (lb) -- (lo);
\draw [->] (lo) -- (mo);
\draw [->] (mo) -- node[above]{$u$} (ro);
\draw [->] (rb) -- (ro);
\draw [->] (mb) -- (mo);
\draw [->] (mo) -- (moo);
\draw [->] (ro) -- (moo);
\draw [->] (lo) -- (loo);
\draw [->] (loo) -- (moo);
\end{tikzpicture}
\end{center}
with each of the three squares a pullback square.

\section{More general fibre bundles}\label{zes}

In this section, we will prove 
a version of theorem \ref{grondsmakreefknol} 
for fibre bundles which are not principal.
It would however be overly optimistic to expect an analogue
of of theorem \ref{grondsmakreefknol} to hold for arbitrary 
smooth fibre bundles, so we will restrict ourselves to those bundles 
that carry a
sufficiently rigid structure on their fibres.




\subsection{Structured fibre bundles}
We start by making this statement more precise.

\begin{Definition} Let $\mathbf{C}$ be a subcategory of the category of smooth manifolds such that 
the group of automorphisms of each object of $\mathbf{C}$  
is a finite dimensional Lie group. 
Then a `structured fibre bundle'
with structure $\mathbf{C}$ and fibre $F_0 \in \mathrm{ob}(\mathbf{C})$ is 
by definition a smooth fibre bundle $\pi : F \rightarrow M$ where the fibres
are objects of $\mathbf{C}$. We also require
each point to possess a neighbourhood $U$ and a local trivialisation 
$\phi :  \pi^{-1}(U) \rightarrow F_0 \times U$ such that $\phi$ restricted to
a single fibre is a $\mathbf{C}$-isomorphism $\pi^{-1}(m) \rightarrow F_0$.  
\end{Definition}

For example, a structured fibre bundle in the category of finite dimensional
vector spaces is a vector bundle.

If $\pi : F \rightarrow M$ is any smooth fibre bundle, then an 
automorphism of $\pi$ 
is by definition a 
diffeomorphism $\alpha$ of $F$ such that  
$\pi(f) = \pi (f')$ implies $\pi(\alpha(f)) = \pi(\alpha(f'))$.
It is called vertical if it maps each fibre to itself.

\begin{Definition} We define an automorphism of a structured fibre bundle to be an automorphism
of the smooth fibre bundle such that its restriction to each single fibre 
is an isomorphism in $\mathbf{C}$.      
\end{Definition}    
    
One can then construct a sequence of groups 
\begin{equation} \label{structuurgroep}
1 \rightarrow \mathrm{Aut}^{\mathbf{C}}_c(F)^{V} 
\rightarrow \mathrm{Aut}^{\mathbf{C}}_c(F)
\rightarrow \mathrm{Diff}_c(M)
\rightarrow 1 
\end{equation}
and its corresponding exact sequence of Lie algebras
\begin{equation} \label{structuurvezel}
0 
\rightarrow \Gamma^{\mathbf{C}}_c(TF)^{V} 
\rightarrow \Gamma^{\mathbf{C}}_c(TF)^{P}
\rightarrow \Gamma_c(TM)
\rightarrow 
0 \,,
\end{equation}
where `$V$' is for vertical, `$P$' for projectable, and $c$
again stands for `$0$ outside a compact subset of $M$'.

The proof of the following corollary of theorem \ref{grondsmakreefknol} 
is now a formality.
\begin{Corollary}\label{structuurcor}
Let $\pi : F \rightarrow M$ be a structured fibre bundle with fibre $F_0$
such that 
(\ref{structuurvezel}) splits as a sequence of Lie algebras.
Then there exists an action $\rho$ of $G(k,M)$ 
by $\mathbf{C}$-automorphisms on a single fibre $F_m$ such that
$$
F = \tf \times_{\rho} F_0 \,.
$$
\end{Corollary}

\proof
Construct the principal $\mathrm{Aut}^{\mathbf{C}}(F_0)$-bundle 
$\pi : P \rightarrow M$, the fibre over $x$
of which is precisely the set of $\mathbf{C}$-isomorphisms 
$\phi : F_{m} \rightarrow F_{x}$.
Then there is a natural isomorphism 
$\mathrm{Aut}^{\mathbf{C}}_{c}(F) \simeq \mathrm{Aut}_{c}(P)^{G}$ under which 
the vertical subgroups of the two correspond, so that the exact sequence 
(\ref{structuurgroep}) is isomorphic to (\ref{flapperdrop}),
and therefore (\ref{structuurvezel}) to (\ref{rijtje}).
As $F = P \times_{\mathrm{Aut}(F_m)} F_{m}$,
we can now apply theorem \ref{grondsmakreefknol} to $P$ in order to substantiate our
claim. \qed


\subsection{Vector bundles}
We specialise to the case of vector bundles. As we pointed out, 
these are precisely structured fibre bundles in the category of finite
dimensional vector spaces.

The exact sequence of Lie algebras (\ref{structuurvezel}) for a vector bundle
$E$ with fibre $V$ is then
\begin{equation}\label{veclie}
1 \rightarrow \mathrm{DO}^{0}_{c}(E) \rightarrow \mathrm{DO}^{1}_{c}(E) 
\rightarrow \mathrm \Gamma_{c}(TM) \rightarrow 0\,,
\end{equation}
where $\mathrm{DO}_{c}^{1}(E)$ is the Lie algebra of compactly supported 
$1^{\mathrm{st}}$ order differential operators on $\Gamma(E)$, and
$\mathrm{DO}_{c}^{0}(E)$ the ideal of $0^{\mathrm{th}}$ order ones, 
that is to say 
$\mathrm{DO}_{c}^{0}(E) \simeq \Gamma_{c} (E \otimes E^{*})$.  

Corollary \ref{structuurcor} then says that (\ref{veclie})
splits as a sequence of Lie algebras if and only if there
is a representation $\rho$ of $G(k,M)$ on $V$ 
such that
$
E \simeq \tf \times_{\rho} V\,.
$

But thanks to the fact that all finite dimensional representations
of the universal cover of $\mathrm{GL}^{+}(\mathbf{R}^{n})$
factor through $\mathrm{GL}^{+}(\mathbf{R}^{n})$ itself, we can 
even say something slightly stronger.

\begin{Proposition} \label{structuurvector}
Let $E \rightarrow M$ be a vector bundle 
for which (\ref{veclie}) splits as a sequence of Lie algebras. 
Then there exists a representation
$\rho$ of $(G^{k} \times \pi_{1}(M))_{m,m}^{\mathrm{Pr}}$ on V such that 
$$
E \simeq (G^{k}(M) \times \pi_{1}(M))_{*,m}^{\mathrm{Pr}} 
\times_{\rho} V\,.
$$
\end{Proposition}
{\bf Remark \quad}
If $M$ is orientable, this reads   
$
E \simeq \pi^{*} F^{+k}(M) 
\times_{\rho} V
$, with $\pi^{*} F^{+k}(M)$ the pullback of $F^{+k}(M)$ along
$\pi : \tilde{M} \rightarrow M$, considered as a principal
$G_{0,0}^{+k}(\mathbf{R}^{n}) \times \pi_{1}(M)$-bundle
over $M$.

\proof
%
Consider the restriction of the map $\tau$ in equation 
(\ref{tauprojectie}) to the group $\tilde{G}_{m,m}^{+k}(M)$.
In order to prove the proposition, we need but show that 
its kernel $Z$ 
acts trivially on $V$.
For $k=0$, this is clear.

If $k$ is at least 1, 
the homomorphism 
$
\tglnp
\rightarrow
\tilde{G}_{m,m}^{+k}(M)
$
makes $V$ into a finite dimensional representation
space for $\tglnp$.  
But it is known (see \cite[p.\,311]{Kn}) that all finite dimensional
representations of its cover factor through $\glnp$
itself. This implies that
the subgroup 
$Z$ which covers the identity
must act trivially on $V$,
and we may consider 
$$\tilde{G}_{*,m}^{+k}(M) / Z 
\simeq 
(G^{k}(M) \times \pi_{1}(M))_{*,m}^{\mathrm{Pr}}
$$
to be the underlying bundle, as announced.\qed

This reduces the problem of classifying vector bundles with split 
sequence (\ref{veclie})
to the representation
theory of $(G^{k} \times \pi_{1}(M))_{m,m}^{\mathrm{Pr}}$.
%

The above extends a result \cite{Te} of Terng, in which she
classifies vector bundles which allow for a local splitting of the sequence of 
groups (\ref{structuurgroep}). 
It is an extension first of all in the sense that 
we prove, rather than assume, that the splitting is local.
Secondly, we have shown that in classifying 
vector bundles
with split sequence (\ref{veclie}) of Lie algebras 
rather than groups, one encounters only slightly more.  
Intuitively speaking, the extra bit is the representation theory
of $\pi_{1}(M)$.
We refer to \cite{Te} for a thorough exposition of the representation theory
of $G^{k}_{0,0}(\mathbf{R}^{n})$.

\section{Flat connections}\label{vijf}

Having concluded our classification of bundles
in which (\ref{rijtje}) is split as a sequence
of Lie algebras, the time has come to
apply our newly acquired knowledge.

In this section, we will investigate splittings
that come from a flat equivariant connection on a 
principal $G$-bundle $P \rightarrow M$.
We will prove that if the Lie algebra 
$\mathfrak{g}$ of $G$
does not contain $\mathfrak{sl}(\mathbf{R}^{n})$
as a subalgebra, 
then 
the sequence of Lie algebras (\ref{rijtje}) splits
if and only if $P$ admits a flat equivariant connection.
In other words, the sequence (\ref{rijtje}) then splits
as a sequence of Lie algebras if and only if it
splits as a sequence of Lie algebras and $C^{\infty}(M)$-modules.

Note that this is certainly not the case for general groups $G$.
The frame bundle for example always allows for a splitting of 
(\ref{rijtje}), but usually not for a flat connection.

\subsection{Lie algebras that do not contain 
$\mathfrak{sl}(\mathbf{R}^{n})$}

Although lemma \ref{hallekal} exhibits $\sigma$
as a differential operator of finite order,  
the bound on the order is certainly not optimal. 
With full knowledge of the Lie algebras
at hand, sharper restrictions can be put 
on the kernel of $\sigma$.
In particular, if $\mathfrak{g}$
does not contain $\mathfrak{sl}(\mathbf{R}^{n})$,
there is only a single relevant ideal, and $\sigma$ is of order 
at most 1.
For notation, see section \ref{localmap}.
%


\begin{Lemma}\label{grondkusrelmuis}
Let $n=1$, and let $\mathfrak{g}$ be such that it does not 
contain two nonzero elements such that $[X,Y] = Y$. 
Or
let $n \geq 2$, and let $\mathfrak{g}$ be such that it does not
admit $\mathfrak{sl}(\mathbf{R}^n)$ as a subalgebra.
Then the kernel of the homomorphism 
$\check{f}_m : \vecn \rightarrow \mathfrak{g}$ 
contains 
$\{v \in \vecn \, | \, \mdiv_m(v) = 0\}$.
\end{Lemma}
\proof
We start with the case $n=1$. Again, we note that the
only ideals of
$\mathrm{Vec}_{1} = \mathrm{Span}\{x^{k} \del_x \,|\, k \geq 1\}$
are 
$\mathrm{Span}\{x^2\del_x , x^k \del_x \,| \, k \geq 4 \}$
and $\mathrm{Span}\{x^k \del_x \,| \, k \geq N \}$ with $N \geq 1$.  
The corresponding quotients all contain $X$ and $Y$ with $[X,Y]=Y$,
except the ones corresponding to the cases $N=1,2$.
This means that also $\mathfrak{g}$, containing
$\mathrm{Vec}_{1} / \ker(\check{f}_m)$ as the image of $\check{f}_m$,  
will possess $X$ and $Y$ such that $[X,Y] = Y$ unless 
the kernel of $\check{f}_m$ contains 
$\mathrm{Span}\{x^k \del_x \,| \, k \geq 2 \} = 
\{v \in \mathrm{Vec}_1 \, | \, \mdiv_m(v) = 0\}$.

Now for $n \geq 2$.
Under the identification $\vecn^{0} \simeq \mathfrak{gl}(\mathbf{R}^n)$
 given by $x_{i}\del_j \mapsto e_{i j}$,
the Euler vector field is the identity $\mathbf{1}$
and $\mdiv_{m}$ becomes the trace. 
As $\ker(\check{f}_m)^{0}$ is an ideal in $\mathfrak{gl}(\mathbf{R}^n)$,
it can be either $0$, $\mathbf{R} \mathbf{1}$, $\mathfrak{sl}(\mathbf{R}^n)$
or $\mathbf{R} \mathbf{1} \oplus \mathfrak{sl}(\mathbf{R}^n)$. 
In the former two cases, 
$\mathrm{Im}(\check{f}_m) \simeq \vecn / \ker(\check{f}_m)$,
and hence $\mathfrak{g}$,  
would contain $\mathfrak{sl}(\mathbf{R}^n)$ as a subalgebra, 
contradicting the hypothesis.
Hence $\mathfrak{sl}(\mathbf{R}^n) \subseteq \ker(\check{f}_m)$.
If we now show that 
$[\vecn , \mathfrak{sl}(\mathbf{R}^n)] = 
\mathfrak{sl}(\mathbf{R}^n) \bigoplus_{k=1}^{\infty} \vecn^{k}$,
the proof will be complete.

Let $i \neq j$. We then have 
$[x_{i}\del_j , x_{j}x^{\vec{\alpha}}\del_{j}] = (\alpha_j + 1) x_i
x^{\vec{\alpha}} \del_j$, showing that 
$x_i x^{\vec{\alpha}} \del_j \in [\vecn , \mathfrak{sl}(\mathbf{R}^n)]$.
The only basis elements not of this shape are of the form $x_{j}^{k}\del_j$.
But $[x_{j}\del_i , x_i x_{j}^{k-1}\del_{j}] = 
x_{j}^{k}\del_j - x_i x_{j}^{k-1}\del_{i}$. If $k \geq 2$, the latter part 
was just shown to be in $[\vecn , \mathfrak{sl}(\mathbf{R}^n)]$, so that also 
$x_{j}^{k}\del_j \in [\vecn , \mathfrak{sl}(\mathbf{R}^n)]$. 
If $k=1$, the elements 
$x_j \del_j - x_i \del_i$ join $x_i \del_j$ to form a basis of 
$\mathfrak{sl}(\mathbf{R}^n)$.
\qed

This rather limits the possibilities. 
Not only can we restrict to first order, but also 
the Lie algebroid map 
$\nabla : J^{1}(TM) \rightarrow TP/G$ vanishes on
the trace-zero jets
$K_{m} = \{\, j_{m}^{1}(v) \in J^{1}(TM) \, | \, v(m) = 0 \, \, \mathrm{and} \,\,
\mathrm{Tr}(v)=0\, \}$,
so that it factors through the `trace Lie algebroid' 
$\mathrm{Tr}_{m}(M) := J_{m}^{1}(TM)/K_{m}$.

This in turn is the Lie algebroid of the `determinant groupoid'
$\mathrm{Det}(M)$. An element $[\alpha]_{m',m}$ of $\mathrm{Det}(M)_{m',m}$
is by definition an equivalence class of diffeomorphisms mapping 
$m$ to $m'$, with $\alpha \sim \beta$ if and only if
$\mathrm{Det}(\beta^{-1}\alpha) = 1$. 

As $[\alpha]_{m',m}$ identifies $\wedge^{n}(T^{*}_{m}M)$ 
with $\wedge^{n}(T^{*}_{m'}M)$,
the source fibre $\mathrm{Det}(M)_{*,m}$ is 
isomorphic to the determinant line bundle $\wedge^{n}(T^{*}M) \rightarrow M$.
Its\footnote{An isomorphism 
$\wedge^{n}(T^{*}M) \simeq \mathrm{Det}(M)_{*,m}$ is only given after 
a choice of $\lambda_{0} \in \wedge^{n}(T_{m}^{*}M)$. This determines 
the connected component.} 
connected component $\wedge^{n,+}(T^{*}M)$ is the 
the bundle of positive top forms
if $M$ is orientable, and the whole bundle otherwise.
 
Its universal covering space is the bundle $\wedge^{n,+}(T^{*}\tilde{M})$
of positive top forms on $\tilde{M}$.
Indeed, 
$\tilde{M}$ is always orientable, regardless of whether or not $M$ is.
This means that 
$\wedge^{n}(T^{*}\tilde{M})$ is a trivial bundle,
and that its connected component 
$\wedge^{n,+}(T^{*}\tilde{M}) \simeq \tilde{M} \times \mathbf{R}^{+}$ 
is simply connected. The covering map is induced by the map 
$\tilde{M} \rightarrow M$. This leads to the following version of 
theorem \ref{grondsmakreefknol}.


\begin{Proposition}\label{habbelsnap}
Let $P$ be a principal $G$-bundle over an $n$-dimensional manifold $M$.
Let $G$ be such that its Lie algebra $\mathfrak{g}$ does not 
contain $\mathfrak{sl}(\mathbf{R}^{n})$ if $n > 1$,
or $[X,Y] = Y$ if $n=1$. 
Then there is a homomorphism $\pi_{1}(M) \times \mathbf{R}^{+} \rightarrow G$
associating $P$ to the principal $\pi_{1}(M) \times \mathbf{R}^{+}$-bundle
$
\wedge^{n,+}(T^{*}\tilde{M}) \rightarrow M 
$.
$$
P \simeq \wedge^{n,+}(T^{*}\tilde{M}) \times_{\pi_{1}(M) \times \mathbf{R}^{+}}G\,.
$$
\end{Proposition}
We may even classify the possible splittings. 
\begin{Corollary}\label{grondsmakgu}
Under the hypotheses of proposition \ref{habbelsnap}, any 
Lie-algebra homomorphism 
$\sigma: \Gamma_{c}(TM) \rightarrow \Gamma_{c}(TP)^{G}$
which splits the sequence of Lie algebras (\ref{rijtje}) 
can be written
\begin{equation}\label{raldaldalgu}
\sigma = \nabla^{\mu} + \Lambda \mdiv_{\mu}\,,
\end{equation}
where $\nabla^{\mu}$ is a flat equivariant connection on $P$,
and $\Lambda$ a section of $\mathrm{ad}(P)$ which is 
constant w.r.t. the connection induced on $\mathrm{ad}(P)$ by $\nabla^{\mu}$. 
\end{Corollary}
{\bf Remark \quad} In particular, this shows that there exists a flat
connection which splits (\ref{rijtje}), even though 
most splittings are not flat connections.

\proof First, we prove the case $P = \wedge^{n,+}(T^{*}\tilde{M})$. 
Pick a nonzero \mbox{(pseudo-)} density $\mu$ on $M$. This induces 
an honest density $\tilde{\mu}$ on $\tilde{M}$, which in turn
identifies 
$
\wedge^{n,+}(T^{*}\tilde{M}) 
$
with
$\tilde{M} \times \mathbf{R}^{+} 
$. 
The local trivialisations of
$\wedge^{n,+}(T^{*}\tilde{M}) \rightarrow \tilde{M}$
and 
$\tilde{M} \rightarrow M$
combine to locally trivialise 
$\wedge^{n,+}(T^{*}\tilde{M}) \rightarrow M$.
This yields a flat 
equivariant
connection $\nabla^{\mu}$ on $\wedge^{n,+}(T^{*}\tilde{M})$,
which annihilates $\tilde{\mu}$.

The splitting $\sigma$ is uniquely determined by the action of $\sigma(v)$ 
on local sections $\tilde{\nu}$, 
which reads 
$\sigma(v) (\tilde{\nu}) = \pi^{*} \circ \mathcal{L}_{v} \circ 
\pi^{* -1}
\tilde{\nu}$,
where $\pi$ is the map from $\tilde{M}$ to $M$. 
If we define the divergence w.r.t. $\mu$
by the requirement that the Lie derivative
$\mathcal{L}_v \mu $ equal $\mathrm{Div}_{\mu}(v) \mu$,
then 
\begin{eqnarray*}
\sigma(v)(f\tilde{\mu}) &=& \pi^{*} (\mathcal{L}_{v}(f\mu))\\
&=& \pi^{*}(v(f) \mu + \mathrm{Div}_{\mu}(v) f \mu)\\
&=& v(f) \tilde{\mu} + \mathrm{Div}_{\mu}(v) f \tilde{\mu}\\
&=& \nabla^{\mu}_{v} (f\tilde{\mu}) + \mathrm{Div}_{\mu}(v) f \tilde{\mu}\,.
\end{eqnarray*}
This shows that $\sigma(v) = \nabla^{\mu}_{v} + \Lambda \mathrm{Div}_{\mu}(v)$, 
with $\Lambda = \del_{r}$, the equivariant vertical vector field defined by the 
action of $\mathbf{R}^{+}$ on $\wedge^{n}(T^{*}\tilde{M})$.
(Equivariant vertical vector fields on $P$ 
correspond to sections of $\mathrm{ad}(P)$.)

The general case follows by proposition \ref{habbelsnap}. \qed

\subsection{Lie algebra cohomology}

If we specialise to the case of a trivial bundle over an 
abelian group $G$, we find ourselves in the realm of Lie algebra cohomology.
The continuous cohomology of the Lie algebra of vector 
fields with values in the functions has already been 
unravelled in all degrees \cite{Fu}.
Corollary \ref{grondsmakgu} describes this cohomology only in
degree 1, but now with all cocycles rather than just the continuous ones.
\begin{Corollary}\label{homologerengu}
Let $H_{LA}$ denote Lie algebra cohomology and $H_{dR}$ de Rham cohomology.
Let $\mathfrak{g}$ be abelian, and consider the representation
$C_{c}^{\infty}(M,\mathfrak{g})$ of $\Gamma_{c}(TM)$ where 
a vector field $v$ acts by the Lie 
derivative $\mathcal{L}_v$. Then 
$$
H^1_{LA}(\Gamma_{c}(TM) ,C_{c}^{\infty}(M,\mathfrak{g}))
\simeq
H^1_{dR}(M,\mathfrak{g}) \oplus \mathfrak{g}\,.
$$
\end{Corollary}
\proof 
Consider the trivial
bundle $M \times G \rightarrow M$ over an abelian Lie group $G$, 
which comes equipped with a 
flat connection $\nabla^{0}$,
which acts as Lie derivative.  
Note that abelian $\mathfrak{g}$ certainly satisfy 
the conditions of propositions \ref{habbelsnap} and \ref{grondsmakgu}. 
View $\Gamma_{c}(\mathrm{ad}(P)) \simeq C^{\infty}_{c}(M,\mathfrak{g})$
as a representation of $\Gamma_{c}(TM)$, and consider its
Lie algebra cohomology. An $n$-cochain is an alternating linear map 
$\Gamma_{c}(TM)^{n} \rightarrow C_{c}^{\infty}(M,\mathfrak{g})$.
For $f^{1} \in C^{1}$, closure $\delta f^{1} = 0$ amounts to
$$
\mathcal{L}_{v}f^{1}(w) - \mathcal{L}_{w}f^{1}(v) - f^{1}([v,w]) 
= 0\,. 
$$
Due to this cocycle condition, $\sigma = \nabla^{0} + f^{1}$ 
is once again a Lie algebra homomorphism splitting $\pi_{*}$.
According to corollary \ref{grondsmakgu},
it must therefore take the shape 
$\sigma = \nabla^{\mu} + \Lambda \mdiv_{\mu}$,
where $\Lambda \in \mathfrak{g}$ is constant.
One can write $\nabla^{\mu} = \nabla^{0} + \omega^{1}$
for some closed 1-form $\omega^{1}$, so that
$f^{1} = \omega^1 + \Lambda \mdiv_{\mu}$.
This classifies the closed $1$-cocycles.


Exact 1-cocycles satisfy $f^1(v) = \delta f^{0} (v) = 
\mathcal{L}_{v}f^{0} = df^{0}(v)$, with $f^0$ a 
$0$-cocycle, that is an element of
$C^{\infty}(M,\mathfrak{g})$.

Note that a change of density $\mu' = e^{h}\mu$
alters $f^1$ by a mere coboundary $\Lambda dh$,
so that the choice of $\mu$ is immaterial. 
The class of $\omega^1 + \Lambda \mdiv_{\mu}$
modulo $\delta C^0$ is therefore determined 
by $[\omega^1] \in H^1_{dR}(M,\mathfrak{g})$
and $\Lambda \in \mathfrak{g}$. 
\qed

Continuity turns out to be implied by the closedness-condition.
A similar situation was encountered by Takens in \cite{Ta}, when proving that
all derivations of $\Gamma_{c}(TM)$ are inner, i.e.
$H^1_{LA}(\Gamma_{c}(TM) ,\Gamma_{c}(TM)) = 0$.

\section{Gauge theory, general relativity and spinors}\label{zeven}
In this section, we briefly reflect on the relationship
between (generalised) spin structures and principal fibre
bundles with a splitting of (\ref{rijtje}). 

We formulate general relativity and gauge theory, 
including the fermionic fields, 
in terms of principal fibre bundles over the 
manifold $M$ which describes space-time.
For convenience, we will take
$M$ to be a smooth and orientable manifold of dimension at least 2.

\subsection{General relativity}

The fundamental degrees of freedom in general relativity
are a pseudo-Rie\-mann\-ian metric $g$ on space-time $M$, 
and a connection $\nabla$ on $TM$. 

We identify the metric
with a section of $F^{+}(M)/\mathrm{SO}(\eta)$ by 
associating to
$g_{x} : T_{x}M \times T_{x}M \rightarrow \mathbf{R}$ the coset 
of all frames $f : \mathbf{R}^{n} \rightarrow T_{x}M$ 
such that $f^{*}g = \eta$.
We identify the connection on $TM$ with an equivariant connection on 
$F^{+}(M)$, that is a section of
$J^{1}(F^{+}(M))/\mathrm{GL}^{+}(\mathbf{R}^{n}) \rightarrow M$.
(Its value at $x$ is the class of a
1-jet at $x$ of a section $\phi \in \Gamma(F^{+})$ with 
$\nabla_{x} \phi = 0$ \cite{Sr}.) 
These two, metric and connection, 
are conveniently combined into
a single section of $J^{1}(F^{+}(M))/\mathrm{SO}(\eta)$. 

The dynamics of the theory are then governed by 
the Einstein-Hilbert action 
$S_{EH} : \Gamma(J^{1}(F^{+}(M))/\mathrm{SO}(\eta)) \rightarrow \mathbf{R}$,
defined in terms of the Ricci scalar $R$ by
$$
S_{EH}(g,\nabla) = \int_{M} R(g,\nabla) \sqrt{|g|}dx_{0}\ldots dx_{n}\,.
$$
In the Einstein-Hilbert approach, the connection is constrained 
to equal the Levi-Civita connection, and only the metric is varied.
In the Palatini approach \cite{As},
the connection varies 
independently, and the fact that $\nabla$ is the Levi-Civita
connection of the metric is a consequence of the field equations.




Either way, it is clear that the fields transform in a 
natural fashion under diffeomorphisms of $M$.
Any diffeomorphism $\alpha$ of $M$ lifts to an 
automorphism $\Sigma(\alpha)$ of 
$J^{1}(F^{+}(M))/\mathrm{SO}(\eta)$, defined
by $\Sigma(\alpha)(j^{1}_{m}(\phi)) = j^{1}_{m}(\alpha_{*}\circ \phi \circ
\alpha^{-1})$.
This splits the sequence of groups
$$
1 
\rightarrow
\mathrm{Aut}^{V}( J^{1}(F^{+}(M))/\mathrm{SO}(\eta)) )
\rightarrow
\mathrm{Aut}(J^{1}(F^{+}(M))/\mathrm{SO}(\eta))
\rightarrow
\mathrm{Diff}(M)
\rightarrow
1\,.
$$

This splitting is central to the theory of general relativity.
The requirement that the action be invariant under co-ordinate 
transformations, $S_{EH}(\Sigma(\alpha) \phi) = S_{EH}(\phi)$,
cannot even be formulated without providing $\Sigma$ explicitly.

Note that as the above sequence of groups splits, so does
the corresponding sequence of Lie algebras of vector fields.





\subsection{Fermions and spin structures}
We wish to describe fermions.
As these are known to transform under Lorentz transformations 
by a projective representation rather than a linear one,
we must extend our framework.
\subsubsection{The bundle}
Suppose that we have a spin structure $Q$ w.r.t
a background metric $g$. 
Let $u : \hat{Q} \rightarrow F^{+}(M)$ be as in section \ref{Spinmanifolds}, and
let $V$ be a unitary spinor representation\footnote{
The indefinite article is appropriate since there is a choice involved here.
The connected component of $\one$ of $\widetilde{\mathrm{SO}}(3,1)$ 
is $\mathrm{Spin}^{\uparrow}(3,1) \simeq \mathrm{SL}^{2}(\mathbf{C})$. 
A spinor representation 
for the connected component can then be unambiguously derived from 
a Clifford algebra representation \cite{He}.
But as $\widetilde{\mathrm{SO}}(3,1)$ is not isomorphic to
$\mathrm{Spin}(3,1)$, the action of the order 2 central elements
covering $PT$ will have to be specified `by hand'.
Again, this is not relevant if $M$ is time-orientable as well as 
orientable.
} 
of 
$\widetilde{\mathrm{SO}}(\eta)$.
Then one has the composite bundle
$$\hat{Q} \times_{\widetilde{\mathrm{SO}}(\eta)} V 
\rightarrow F^{+}(M)/\mathrm{SO}(\eta)
\rightarrow M\,.
$$
A section $\tau : M \rightarrow \hat{Q} \times_{\widetilde{\mathrm{SO}}(\eta)} V$
can then be interpreted as a metric $g$
along with a fermionic field $\psi$.
Consider $g$ as the induced section of
$F^{+}(M)/\mathrm{SO}(\eta)$.  
Use $g$ to construct the spinor bundle
$u^{-1}(OF^{+}_{g}) 
\times_{\widetilde{\mathrm{SO}}(\eta)} V$,
and obtain a section $\psi$  
by simply restricting the image of $\tau$.
 
In the same vein, we will take a physical field to be a section of
\mbox{$J^{1}(\hat{Q})\times_{\widetilde{\mathrm{SO}}(\eta)}V$}. 
This is equivalent to providing a triple of sections, one of
$F^{+}(M)/\mathrm{SO}(\eta)$, one of
$J^{1}(F^{+}(M))/\mathrm{GL}^{+}(\mathbf{R}^{n})$
and one of
$u^{-1}(OF^{+}_{g}) \times_{\widetilde{\mathrm{SO}}(\eta)} V$.
These correspond to  
the metric $g_{\mu \nu}$,
the Levi Civita-connection $\Gamma^{\alpha}_{\mu \beta}$,
and the 
fermionic field $\psi^{i}$ respectively. 
\subsubsection{Transformation behaviour}
Let us investigate its transformation behaviour.
As a spinor changes sign under a $2\pi$-rotation, there is no hope
of finding an interesting homomorphism of groups 
$\mathrm{Diff}(M) \rightarrow \mathrm{Aut}
(J^{1}(\hat{Q})\times_{\widetilde{\mathrm{SO}}(\eta)}V)$.
There is however a canonical homomorphism of Lie algebras.


Because $u : \hat{Q} \rightarrow F^{+}(M)$ has discrete fibres, it
has a unique flat equivariant connection $\nabla$.
This means that the exact sequence of Lie algebras
$$
0\rightarrow
\Gamma_{c}(\mathrm{ad}(\hat{Q}))
\rightarrow
\Gamma_{c}(T\hat{Q})^{\tglnp}
\rightarrow
\Gamma_{c}(TM)
\rightarrow 0
$$
is split by $\sigma := \nabla \circ D$, with
$D : \Gamma_{c}(TM) \rightarrow \Gamma_{c}
(TF^{+}(M))^{\glnp}$
the first order derivative $D(v)_f = \del/\del t|_{0} \exp(tv)_{*}\circ f$.
This induces a splitting for $J^{1}(\hat{Q})$ by prolongation 
(see e.g. \cite{FR}), and consequently also one for
$J^{1}(\hat{Q})\times_{\widetilde{\mathrm{SO}}(\eta)}V$.
Note that $u_{*}\circ \sigma $ equals $D$.

We see that $\hat{Q}$ is an infinitesimally natural bundle,
and that the canonical splitting $\sigma$ of (\ref{rijtje})  
does not come from a splitting of groups.
But even if a (different) splitting does exist at the level of groups, 
it will be irrelevant to physics. 
Take for example
the trivial spin structure 
$Q = \mathbf{R}^{n} \times \widetilde{\mathrm{SO}}(\eta)$ 
over $\mathbf{R}^{n}$,
and 
lift $\alpha \in \mathrm{Diff}(M)$ 
to $\mathrm{Aut}(Q)$ by $\Sigma(\alpha)(m,q) = (\alpha(m),q)$.
This is clearly the wrong thing to do. First of all, sections of
$Q \times_{\widetilde{\mathrm{SO}}(\eta)}V$ transform under the trivial
representation of the Lorentz group. Our fermions are Lorentz scalars
rather than spin-$1/2$ particles.
Secondly, using the wrong splitting will generally
result in an incorrect energy-momentum tensor \cite{GM}.

We arrive at the conclusion that not only the bundle $Q$ and
the covering map $u : Q \rightarrow OF_{g}$ are physically 
relevant, but 
also the splitting 
$\sigma : \Gamma_{c}(TM) \rightarrow \Gamma_{c}(T\hat{Q})^{\tglnp}$.
It must satisfy $u_{*} \circ \sigma = D$
in order for the metric $g \in \Gamma(TF^{+}(M)/SO(\eta))$ to transform
properly.
Such a $\sigma$ is naturally associated to any ordinary
spin structure $Q$. For generalised spin structures however,
this is no longer so.

\subsection{Gauge fields and $\spg$-structures}

In the presence of gauge fields, the topological conditions 
on $M$ in order to support a spin structure are more relaxed.
Intuitively, this is because the gauge group $G$ can absorb some of the 
indeterminacy that stems from the 2:1 cover of the Lorentz group.

\subsubsection{$\spg$-structures}
This is made more rigorous by the notion of a $\spg$-structure \cite{AI}.

\begin{Definition} Let $G$ be a Lie group with a central subgroup isomorphic to $Z$. 
Let\footnote{This notation is convenient but slightly misleading.
Be ware that 
if $\eta$ is of signature $+++-$, then
$\mathrm{Spin}^{Z}$ is isomorphic to 
$\widetilde{\mathrm{SO}}(\eta)$, not to $\mathrm{Spin}(3,1)$. 
} $\spg := \widetilde{\mathrm{SO}}(\eta) \times_{Z} G$.
A $\spg$-structure is a
$\spg$-bundle $Q$ over $M$,
together with a map  $u : Q \rightarrow OF^{+}_{g}$
that makes 
\begin{center}
\begin{tikzpicture}
\pgfsetzvec{\pgfpoint{0.385cm}{-0.385cm}}
\node (lg) at (-1.2,0.7) {$\spg$};
\node (rg) at (1.2,0.7) {$\mathrm{SO}(\eta)$};
\node  (rpijl) at (1.2,0.3) {$\curvearrowleft$};
\node  (lpijl) at (-1.2,0.3) {$\curvearrowleft$};
\node (lb) at (-1.2,0) {$Q$};
\node (rb) at (1.2,0) {$OF^{+}_{g}$};
\node (lo) at (0,-1.2) {$M$};
\draw [->] (lg) --node[above]{$\kappa$} (rg);
\draw [->] (lb) --node[above]{$u$} (rb);
\draw [->] (lb) -- (lo);
\draw [->] (rb) -- (lo);
\end{tikzpicture}
\end{center}
commute. We again denote the map $(x , g) \mapsto \kappa(x)$ by $\kappa$.    
\end{Definition}

This of course gives rise to 
the principal 
$\widetilde{\mathrm{GL}}(\mathbf{R}^{n}) \times_{Z} G$-bundle
$\hat{Q} := Q \times_{\widetilde{\mathrm{SO}}(\eta)}
\widetilde{\mathrm{GL}}(\mathbf{R}^{n})
$. 
If $G = Z$, we recover the notion of a spin structure.
A $\spg$-bundle for the group $U(1)$ is usually called
a spin$^{c}$-structure.

Let $V$ be a representation of $\spg$.
The `physical bundle' is then 
$J^{1}(\hat{Q})\times_{\spg}V$.
A single section of $J^{1}(\hat{Q}) \times_{\spg}V$
represents  
a metric,
a Levi Civita-connection,
a gauge field  
and a fermionic field.

The metric is the induced section of
$F^{+}(M)/\mathrm{SO}(\eta)$,
and the the Levi-Civita connection 
that of
$J^{1}(F^{+}(M))/\mathrm{GL}^{+}(\mathbf{R}^{n})$.
One constructs the 
principal $G/Z$-bundle 
$P := \hat{Q} / \tglnp$, and the gauge field
is the induced equivariant connection on $P$,
a section of $J^{1}(P)/(G/Z)$.
The fermionic field is the induced section of
$\pi^{-1}(OF^{+}_{g}) 
\times_{\widetilde{\mathrm{SO}}(\eta) \times_{Z}G} V$,
where one should note that the bundle
itself depends on $g$.

\subsubsection{Infinitesimally natural $\spg$-structures}

We argue that it only makes sense to consider $\spg$-structures 
which admit an appropriate transformation
law under infinitesimal diffeomorphisms of
space-time. We will call these $\spg$-structures
infinitesimally natural.

\begin{Definition}
A $\spg$ structure $u : Q \rightarrow FO_{g}$ will be called
`infinitesimally natural'
if 
$\hat{Q} \rightarrow M$ is infinitesimally natural as a principal fibre
bundle. Moreover, we require that
the splitting 
$\sigma : \Gamma_{c}(TM) \rightarrow 
\Gamma_{c}(T\hat{Q})^{\widetilde{\mathrm{GL}}(\mathbf{R}^{n})\times_{Z}G}$ 
of (\ref{rijtje}) satisfy $u_{*} \circ \sigma = D$.
\end{Definition}

The mathematical requirement that 
$\hat{Q} \rightarrow M$ be infinitesimally natural as a principal fibre
bundle
corresponds to the physical requirement that
fields should have a definite transformation behaviour under 
infinitesimal co-ordinate transformations.

The requirement $u_{*} \circ \sigma = D$ 
corresponds to the fact that we need to interpret a section of
$\hat{Q}/\widetilde{\mathrm{GL}}(\mathbf{R}^{n})\times_{Z}G 
\simeq F^{+}(M)/\mathrm{SO}(\eta)$ as a metric, and we know that 
its transformation behaviour is governed by $D$.

We view infinitesimally natural $\spg$-structures as the underlying principal
fibre bundles in any classical field theory combining gravity, 
fermions and gauge fields.
Let us work towards 
their classification.

\subsubsection{Classification of infinitesimally natural $\spg$-structures}
Theorem \ref{grondsmakreefknol} is of course the main tool
when classifying infinitesimally natural $\spg$-structures. 
It provides a homomorphism
$
\rho :
G(k,M)
\rightarrow
\widetilde{\mathrm{GL}}(\mathbf{R}^{n})\times_{Z}G
$
such that
$\hat{Q}$ 
is isomorphic to
$
\tf
\times_{\rho} 
(\widetilde{\mathrm{GL}}(\mathbf{R}^{n})\times_{Z}G)
$, and a map $\exp\nabla : \tf \rightarrow \hat{Q}$.
The splitting $\sigma$ is induced by the lift of the $k^{\mathrm{th}}$
order derivative
$
\tilde{D} :
\Gamma(TM) \rightarrow 
\Gamma(T \tf) ^{G(k,M)}
$.
In summary, we have the following commutative diagram.
\begin{center} 
\begin{tikzpicture}
\pgfsetzvec{\pgfpoint{0.385cm}{-0.385cm}}
\node (linkser) at (0, 0) {$0$};
\node (links) at (2.3, 0) 
{$\Gamma_{c}(T\hat{Q})_{v}^{\widetilde{\mathrm{GL}}(\mathbf{R}^{n})\times_{Z}G}$};
\node (midden) at (6, 0) 
{$\Gamma_{c}(T\hat{Q})^{\widetilde{\mathrm{GL}}(\mathbf{R}^{n})\times_{Z}G}$};
\node (rechts) at (9.5,0) {$\Gamma_{c}(TM)$};
\node (rechtser) at (11, 0) {$0$};
\draw [->] (linkser) -- (links);
\draw [->] (links) -- (midden);
\draw [->] (midden) -- node[above] {$\pi_{*}$} (rechts);
\draw [->] (rechts)..controls(8,0.7)..node[above]{$\sigma$}(midden); 
\draw [->] (rechts) --  (rechtser);
\node (2linkser) at (0, -1.2) {$0$};
\node (2links) at (2.3, -1.2) 
{$\Gamma_{c}(TF(M))_{v}^{\mathrm{GL}(\mathbf{R}^{n})}$};
\node (2midden) at (6, -1.2) 
{$\Gamma_{c}(TF(M))^{\mathrm{GL}(\mathbf{R}^{n})}$};
\node (2rechts) at (9.5,-1.2) {$\Gamma_{c}(TM)$};
\node (2rechtser) at (11, -1.2) {$0$};
\draw [->] (2linkser) -- (2links);
\draw [->] (2links) -- (2midden);
\draw [->] (2midden) -- node[above] {$\pi_{*}$} (2rechts);
\draw [->] (2rechts)..controls(8,-0.5)..node[above]{$D$}(2midden); 
\draw [->] (2rechts) --  (2rechtser);
\node (3linkser) at (0, 1.2) {$0$};
\node (3links) at (2.3, 1.2) 
{$\Gamma_{c}(T\tf)_{v}^{G(k,m)}$};
\node (3midden) at (6, 1.2) 
{$\Gamma_{c}(T\tf)^{G(k,m)}$};
\node (3rechts) at (9.5,1.2) {$\Gamma_{c}(TM)$};
\node (3rechtser) at (11, 1.2) {$0$};
\draw [->] (3linkser) -- (3links);
\draw [->] (3links) -- (3midden);
\draw [->] (3midden) -- node[above] {$\pi_{*}$} (3rechts);
\draw [->] (3rechts)..controls(8,1.9)..node[above]{$\tilde{D}$}(3midden); 
\draw [->] (3rechts) --  (3rechtser);
\draw [->] (links) --node[left]{$u_{*}$} (2links);
\draw [->] (midden) --node[left]{$u_{*}$} (2midden);
\draw [->] (rechts) --node[right]{id} (2rechts);
\draw [->] (3links) --node[left]{$\exp\nabla_{*}$} (links);
\draw [->] (3midden) --node[left]{$\exp\nabla_{*}$} (midden);
\draw [->] (3rechts) --node[right]{id} (rechts);
\end{tikzpicture}
\end{center}

The classification theorem for  
infinitesimally natural $\spg$-structures
will take the following form. 
\begin{Theorem}\label{naturalspin}
Let $M$ be an orientable smooth manifold of dimension $n \geq 2$. 
Let $G$ be a Lie group with a central 
subgroup $Z$ isomorphic to $\pi_{1}(\gln)$, and let $\mathrm{Lie}(G)$
be such that it does
not contain $\mathfrak{sl}(\mathbf{R}^{n})$. 
Finally, let $(Q,u)$ be an infinitesimally natural $\spg$-structure
over $M$.
Then $i_{*} : Z \rightarrow \pi_{1}(F^{+}(M))$
is injective, and
there exists a homomorphism
$
\tau : \pi_{1}(F^{+}(M)) \rightarrow G
$
which fixes $Z$, and makes $(Q,u)$ isomorphic to the bundle
$$
\widetilde{OF}^{+}_{g} \times_{\tau} G\,,
$$
with $u$ the natural projection map onto $OF^{+}_{g}$.
\end{Theorem} 
{\bf Remark\quad} The restriction that $\mathrm{Lie}(G)$ not contain 
$\mathfrak{sl}(\mathbf{R}^{n})$ does not appear to be very limiting
as far as gauge groups are concerned. Indeed, one would expect 
a gauge group
to possess a faithful finite dimensional unitary representation
and therefore be compact.
This cannot be if $\mathfrak{sl}(\mathbf{R}^{n})$ is a subalgebra
of $\mathrm{Lie}(G)$.

{\bf Remark\quad} We have already seen that an ordinary spin structure is always 
infinitesimally natural, so that proposition \ref{klopboor?} is just the special case $G = Z$.

\subsubsection{Two lemmas}

We set out to prove theorem \ref{naturalspin}. We start with 
two lemmas designed to explicate the homomorphism
$
\rho :
G(k,M)
\rightarrow
\widetilde{\mathrm{GL}}(\mathbf{R}^{n})\times_{Z}G
$.


\begin{Lemma}
Denote the natural projection map
$\tf \rightarrow F^{+}(M)$ by $\pi$,
and write $\nu : G(k,M) \rightarrow \glnp$ for the corresponding homomorphism of groups.
Then each infinitesimally natural $\spg$-structure is isomorphic to one
for which 
$\kappa \circ \rho = \nu$
and
$u \circ \exp\nabla = \pi$.
\end{Lemma}
\proof 
Denote $u \circ \exp\nabla$ by $\gamma$. 
Pick $\tilde{f}  \in \tfx$,
and let $\pi(\tilde{f}) = f$.
There exists a $c \in \glnp$
such that $\gamma(\tilde{f}) = f c$.
Equivariance of $\gamma$ then implies
$\gamma(\tilde{f} \tilde{y}) = f c \kappa(\rho (\tilde{y}))$
for all $\tilde{y} \in G(k,M)$.

We required $\gamma_{*} \circ \tilde{D} = D$, but 
clearly we also have $\pi_{*} \circ \tilde{D} = D$. 
Therefore,
if $\exp(\tilde{D}(v)) \tilde{f} = \tilde{f} \tilde{y}$
and 
$\exp({D}(v)) f = f y$,
we must have both 
$\nu(\tilde{y}) = y$
and
$
\kappa(\rho(\tilde{y})) = c^{-1} y c
$.

We may obtain all $\tilde{y}$ in $G(k,M)^{0}$,
the connected component of $G(k,M)$,
by choosing an appropriate $v$, so that
we have 
$
\kappa(\rho(\tilde{y})) = c^{-1} \nu(\tilde{y}) c
$
for all 
$\tilde{y} \in G(k,M)^{0}$.
Although $c$ depends on $\tilde{f}$ a priori,
it turns out to be constant up to scaling.
Indeed, as both $\kappa$ and $\rho$ are constant, 
so is the adjoint action of $c$.
Moreover, $\tilde{f} \mapsto c(\tilde{f})$ is invariant under
$G(k,M)^{0}$,
making it a function on 
$\tf / G(k,M)^{0} \simeq
\tilde{M}$.
All in all, we have established that there exist
$c_{0} \in \mathrm{SL}(\mathbf{R}^{n})$
and
$h : \tilde{M} \rightarrow \mathbf{R}^{+}$
such that
$c(\tilde{f}) = c_{0} h(\tilde{x})$, with $\tilde{x}$
the projection of $\tilde{f}$ to $\tilde{M}$.
We may write 
$u \circ \exp\nabla(\tilde{f}) = \pi(\tilde{f})c_{0}h(\tilde{x})$.


We show that we may as well take $c$ and $h$ to be $1$. 
Pick a $\tilde{c} \in \widetilde{\mathrm{SL}}(\mathbf{R}^{n})$ 
which covers $c$, and
construct the bundle $\hat{Q}_{c} := 
\hat{Q} \times_{\mathrm{Ad}(\tilde{c})}
\widetilde{\mathrm{GL}}(\mathbf{R}^{n})$.
It is isomorphic to $\hat{Q}$,
with isomorphism $\hat{Q}_{c} \rightarrow \hat{Q}$ given by 
$[q,y] \mapsto q\tilde{c}^{-1}h^{-1}(\tilde{x})y$.

If we simply pull back the covering map on $Q$, we obtain
$u_c : \hat{Q}_{c} \rightarrow F^{+}(M)$
given by $[q,y] \mapsto u(q\tilde{c}^{-1}h^{-1}(\tilde{x})y)$.
This makes $Q_{c} := u_{c}^{-1}(OF_{g})$
into a $\spg$-structure 
isomorphic to $Q$, but with the desired properties.
\qed

Recall from (\ref{orient}) that 
$G(k,M) \simeq G(k,\mathbf{R}^{n}) \times_{Z} \pi_{1}(F^{+}(M))$,
and that $G(k,\mathbf{R}^{n}) \simeq \tglnp \ltimes G^{>1}$,
with $G^{>1}$ the subgroup of $k$-jets 
that are the identity to first order.
We unravel $\rho$, considering it as a map
$$
\tsln \times \mathbf{R}^{+} \ltimes G^{>1} 
\times_{Z} \pi_{1}(F^{+}(M))
\rightarrow
\tsln \times \mathbf{R}^{+} \times_{Z} G\,.
$$

\begin{Lemma}
Under the assumptions of theorem \ref{naturalspin}, the map
$\rho$ 
is completely determined by a homomorphism
$\bar{\rho}_{1} : \pi_{1}(F^{+}(M)) \rightarrow  G$ fixing $Z$,
a homomorphism $\bar{\rho}_{2} : \pi_{1}(F^{+}(M)) \rightarrow \mathbf{R}^{+}$,
and an element $\Lambda$ of $\mathrm{Lie}(G)$ which
commutes with $\mathrm{Im}(\bar{\rho}_{1})$. 
We have
$$
\rho : \,\,(\tilde{x} , e^{t} ,g, [p]) \mapsto 
(\tilde{x} , e^{t} \bar{\rho}_{2}([p]), e^{t\Lambda} \bar{\rho}_{1}([p]))\,.
$$
\end{Lemma}
\proof
Consider $\dot{\rho}$ as a Lie algebra homomorphism 
from
$\lsln \times \mathbf{R} \ltimes \mathfrak{g}^{>1}$
to
$\lsln \times \mathbf{R} \times \mathrm{Lie}(G)$,
and let $\dot{\rho}_{ij}$ be its $(i,j)$
component for $i,j \in \{1,2,3\}$. 

Because $\lsln$ is simple and not contained in $\mathrm{Lie}(G)$,
we must have $\dot{\rho}_{13} = 0$.
Due to the previous lemma, $\dot{\rho}_{12} = 0$ and 
$\dot{\rho}_{11} = \mathrm{id}$. As $\tsln$ is simply connected, 
we must have $(\tilde{x},1,1,1) \mapsto (\tilde{x},1,1)$.
In particular, this forces $i_{*} : Z \rightarrow \pi_{1}(F^{+}(M))$
to be injective.

Again due to the previous lemma, $\dot{\rho}_{21} = 0$
and $\dot{\rho}_{22} = \mathrm{id}$. 
Define the `scaling element' $\Lambda := \dot{\rho}_{23}(1)$.
Then $(1,e^{t},1,1) \mapsto (1,e^{t},e^{t\Lambda})$.
The image of $\pi_{1}(F^{+}(M))$ must commute with 
$(\tilde{x} , e^{t} , e^{t\Lambda})$, so that 
$\bar{\rho} : \pi_{1}(F^{+}(M)) \rightarrow \mathbf{R^{+}} \times G$
is well defined, and $\Lambda$ commutes with its image.

We now show that 
$\dot{\rho}(\mathfrak{g}^{> 1}) = 0$.
First of all, as $[\lsln , \lsln + \mathfrak{g}^{> 1}]$
equals
$\lsln + \mathfrak{g}^{> 1}$ if $n$ is at least 2 
(see Lemma \ref{grondkusrelmuis}),
we must have
$$
\dot{\rho}(\mathfrak{g}^{>1}) \subset 
[\dot{\rho}(\lsln) , \dot{\rho}(\lsln + \mathfrak{g}^{>1})]
\subset \lsln \oplus 0 \oplus 0\,.
$$
But on the other hand, as 
$[\mathbf{R} , \mathfrak{g}^{>1}] = \mathfrak{g}^{>1}$,
(recall that $\mathbf{R}$ represents the Euler vector field),
we have 
$$
\dot{\rho}(\mathfrak{g}^{>1}) = 
[\dot{\rho}(\mathbf{R}), \dot{\rho}(\mathfrak{g}^{>1})]
\subset 0 \oplus 0 \oplus \mathrm{Lie}(G)\,.
$$
The intersection being zero, we have $\dot{\rho}(\mathfrak{g}^{>1}) = 0$.
But then $\rho(G^{>1}) = 1$, because $G^{>1}$ is simply connected.\qed

\subsubsection{Proof of theorem \ref{naturalspin}}

Clearly, $\mathbf{R}^{+}$ cannot cause any topological obstruction, so we should
be able to eliminate both $\bar{\rho}_{2}$ and $\Lambda$ from the story.
The former is easy. recall that we seek $Q$, not $\hat{Q}$. 
As $Q$ is a subbundle of $\hat{Q}/\mathbf{R}^{+}$, we will focus on
the latter from now on, allowing us to simply disregard $\bar{\rho}_{2}$.

In order to remove $\Lambda$, we choose a 
volume form $\lambda$ on $M$.
This endows each frame $f \in F_{x}(M)$
with a volume $\mathrm{vol}_{\lambda}(f)$. 
A frame has volume 1 precisely when it is the jet of a diffeomorphism 
which preserves $\lambda$. Denote by $F^{\lambda}(M)$
the $\mathrm{SL}(\mathbf{R}^{n})$-bundle of frames with volume 1.
Its universal cover $\widetilde{F}^{\lambda}(M)$ 
is its inverse image under $\pi$, a principal
$\widetilde{\mathrm{SL}}(\mathbf{R}^{n}) \times_{Z} \pi_{1}(F^{+}(M))$-bundle.

%
%

Define the isomorphism 
$$
\hat{Q}
\rightarrow 
\widetilde{F}_{\lambda}(M) \times_{\pi_{1}(F^{+}(M))} (G \times \mathbf{R}^{+})
$$
by 
$$
(\tilde{f} , g) \mapsto 
\Big(\tilde{f} \, \mathrm{vol}^{-1}_{\lambda}(\pi(\tilde{f})) \,,\,
\exp\big(\log(\mathrm{vol}_{\lambda}(\pi(\tilde{f}))) \,\Lambda \big) g 
, \mathrm{vol}_{\lambda}(\pi(\tilde{f}))\Big)
$$
where we consider $\hat{Q}$ as 
$(\tf \times_{\bar{\rho}_1 \bar{\rho}_2} G)$.
One can see that it is well defined, and that it intertwines the 
natural maps to $F/\mathbf{R}^{+}$.

This shows that $\hat{Q}/\mathbf{R}^{+}$, and therefore the spin structure 
$(Q,u)$, is completely determined 
by the homomorphism $\bar{\rho}_{1} : \pi_{1}(F^{+}(M)) \rightarrow G$.
We denote it by $\tau$ from now on.

If we choose $\lambda$ to be the volume form induced by the metric $g$,
then $OF^{+}_{g}$ is a subbundle of $F_{\lambda}(M)$. Since 
$\pi_{1}(F^{+}(M)) = \pi_{1}(OF_{g}^{+})$, we have
$Q \simeq \widetilde{OF}_{g}^{+} \times_{\tau} G$.
The spin map $u$ is simply the projection 
$\widetilde{OF}^{+}_{g} \rightarrow OF^{+}_{g}$.
The principal $G$-bundle $Q \rightarrow OF_{g}$ has a flat equivariant
connection induced by the one on $\widetilde{OF}^{+}_{g} \rightarrow OF^{+}_{g}$.
This concludes the proof of theorem \ref{naturalspin}. \qed

Tracking back through the isomorphisms, 
we can formulate the following.
\begin{Corollary}
Under the assumptions of theorem \ref{naturalspin},
there exists a $\Lambda \in \mathrm{Lie}(G)$
which com\-mutes with the image of $\tau$,
such that the splitting
$\sigma : \Gamma(TM) \rightarrow
\Gamma(T\hat{Q}/\mathbf{R}^{+})^{\widetilde{\mathrm{SL}}
(\mathbf{R}^{n})\times_{Z}G}$
is given by 
$$
\sigma(v) = \nabla \circ \bar{D}(v) + \mathrm{Div}_{\lambda}(v) \Lambda\,,
$$
where $\lambda$ is the volume form induced by $g$,
$\bar{D}$ is the natural lift
from $\Gamma(TM)$ to $\Gamma(\widetilde{F}^{+}(M))/\mathbf{R}^{+})$,
and we have identified $\Lambda$ with the vector 
field on $\hat{Q}/\mathbf{R}^{+}$
induced by the action of the Lie algebra element.
\end{Corollary}
It is clear that two different homomorphisms 
$\tau_1$ and $\tau_{2}$ yield isomorphic
$\spg$-structures if one is obtained from the other by conjugation
within $G$.

\subsection{Some physical theories} 
If we accept that any physically relevant 
$\spg$-structure must be
infinitesimally natural, and that $G$ should be 
the gauge group of the theory, then
theorem \ref{naturalspin}
provides a link between
the spectrum of elementary particles  
and the global topology of space-time.

The fact that the mere existence of a generalised spin structure
may place restrictions on the space-time manifold was recognised 
by Hawking and Pope \cite{HP}.   
Generalised spin structures were classified \cite{AI},
and it was found that 
if the Lie group contains $\mathrm{SU}(2)$, 
then `universal spin structures' exist 
\cite{BFF}, irrespective of the topology of $M$.
In particular, there are no topological obstructions to the existence of
a $\spg$-structure as soon as $\mathrm{SU}(2) < G$.

But according to theorem \ref{naturalspin}, this changes 
if one requires the $\spg$-structure to be infinitesimally natural. 
Universal spin structures then exist only for certain noncompact 
groups.
For compact $G$, the requirement that there exist a 
homomorphism $\pi_{1}(F^{+}(M)) \rightarrow G$ fixing $Z$
provides an obstruction on the space-time manifold $M$
in terms of the group of internal symmetries $G$. 
Let us see what this means in some specific cases.

\subsubsection{Weyl spinors}
Consider a single massless charged Weyl spinor coupled to 
a $U(1)$ gauge field.
For simplicity, let us assume that $M$ is 4-dimensional, 
oriented, and time-oriented, so that $Z = \mathbf{Z}/2\mathbf{Z}$
and we may use $\mathrm{SL}(\mathbf{C}^{2})$ instead of 
$\widetilde{\mathrm{SO}}(\eta)$. 

This means that $G = U(1)$ and $V = \mathbf{C}^{2} \otimes \mathbf{C}_{q}$,
the two-dimensional
defining representation of 
$\mathrm{SL}(\mathbf{C}^{2})$ tensored with
the one dimensional defining representation 
of $U(1)$. 

Note that the 
representation descends to
$\mathrm{Spin}^{c} = \mathrm{SL}(\mathbf{C}^{2}) \times_{Z} U(1)$.
($Z$ is just $\pm\one$ in $U(1)$.)
This means that upon choosing a $\mathrm{Spin}^{c}$-structure $Q$,
the configuration space is
$\Gamma(J^{1}(\hat{Q})\times_{ \mathrm{ Spin }^{c} }V)$.

%
%
%

Let us now impose the requirement that $Q$ be infinitesimally natural.
Theorem \ref{naturalspin} tells us that there must then be 
a homomorphism $\pi_{1}(F^{+}(M)) \rightarrow U(1)$
sending $\mathbf{Z}/2\mathbf{Z}$ to $\pm 1$. 
If $\pi_{1}(M)$ is finitely generated,
then the image of $\pi_{1}(F^{+}(M))$ in $U(1)$ is a finitely generated
subgroup
containing $\pm\one$. 
It must be isomorphic to $\mathbf{Z}^{n} \times (\mathbf{Z}/2m\mathbf{Z})$
for some $n, m \in \mathbf{N}$. 
If we then send $\mathbf{Z}^{n}$ to $1$, we obtain a homomorphism
$\pi_{1}(F^{+}(M)) \rightarrow \mathbf{Z}/2m\mathbf{Z}$.
The sequence
$1 \rightarrow
\mathbf{Z} / 2\mathbf{Z} \rightarrow
\mathbf{Z}/2m\mathbf{Z} \rightarrow
\mathbf{Z}/m\mathbf{Z}
\rightarrow 1$
splits precisely when $m$ is odd.

We conclude that an infinitesimally natural $\mathrm{Spin}^{c}$-structure
exists on $M$ if and only if there is a surjective homomorphism
$\pi_{1}(F^{+}(M)) \rightarrow \mathbf{Z}/2m\mathbf{Z}$
which preserves $\mathbf{Z}/2\mathbf{Z}$. 
Only if $m$ is odd does this give rise to a spin
structure.

\subsubsection{Dirac spinors}
Next, consider the case of a Dirac spinor. 
That is, $V = \mathbf{C}^{4} \otimes \mathbf{C}_{q}$,
with $\mathbf{C}^{4}$ the representation 
of $\mathrm{Cl}(\mathbf{C}^{4})$ which splits into
two identical irreps
$\mathbf{C}^{2} \oplus \mathbf{C}^{2}$
under $\mathrm{SL}(\mathbf{C}^{2})$, the left handed
and right handed spinors.

The fact that $V$ is reducible under $\mathrm{Spin}^{c}$
makes us re-examine our assumption that
the group $G$ in theorem \ref{naturalspin}
should be the gauge group $U(1)$.
Indeed, the unitary commutant of $\mathrm{SL}(\mathbf{C}^{2})$
in $V$ is $U(2)$ rather than $U(1)$. 
If we take any discrete subgroup $H < U(2)$ and form
the group $U(1)_{H}$ generated by $H$ and $U(1)$,  
can we take $Q$ to
be a $\spg$-structure with structure group $U(1)_{H}$?

As far as only the kinematics is concerned, the answer is yes.
The generic fibre of 
$J^{1}(\hat{Q}) \times_{\spg} V$ is the same for 
$G = U(1)$ as it is for $G = U(1)_{H}$, so 
adding $H$ will not change the space of local sections.

But if we take into account the dynamics, the answer becomes:
`only if $H$ is a group of symmetries of the Lagrangian'.
The reason is that Lagrangians are usually defined 
in local co-ordinates,
yielding a local action $S_{U}$ for each co-ordinate patch $U \subset M$.
However, as the dynamics of the theory should be governed by a global 
action functional $S_{M}$,
it is necessary for $S_{U}$ and $S_{V}$ to agree on $U \cap V$.
This means that if $H$ is part of the structure group of the bundle,
then it must leave the Lagrangian invariant. 
If $H$ is a global symmetry, then the transition functions
must be constant. This is automatic if $H$ is discrete.

For instance, in the case of 
a massive Dirac fermion, the 
subgroup of $U(2)$ which preserves the Lagrangian is precisely the diagonal
$U(1)$.
This means that the relevant 
$\spg$-structures are precisely the $\mathrm{Spin}^{c}$-structures 
classified above.

For massless Dirac spinors, the left and right Weyl spinors decouple,
so that the relevant symmetry group is $U_{L}(1) \times U_{R}(1)$. 
Although the requirement on a manifold to carry 
a $\spg$-structure does not change, this does give us more 
$\spg$-structures for the same manifold.

This illustrates that we may enlarge the gauge group $G$
by any group of discrete symmetries of the Lagrangian in order
to obtain $\spg$-structures.
In particular, this means that infinitesimally natural $\mathrm{Spin}^{c}$-structures
are allowed even for uncharged Weyl spinors. (The image of 
$\pi_{1}(F^{+}(M))$ is automatically discrete.)

\subsubsection{The standard model}
In the the standard model
of elementary particle physics,
the gauge group $G$ is 
$(\mathrm{SU}(3) \times \mathrm{SU}(2)_{L} \times \mathrm{U}(1)_{Y})/ N$,
with $N$ the cyclic subgroup of order 6 generated by
$(e^{2\pi i/3}\one ,  -\one , e^{2\pi i/6})$. 
It is isomorphic to $S(U(3) \times U(2))$, a subgroup of $\mathrm{SU}(5)$,
and it has a unique central subgroup of order 2 
generated by $\mathrm{\textsl{diag}}\,(1,1,1,-1,-1)$.


The fermion representation $V$ for a single family 
can be conveniently described by
$\mathbf{C}^{2} \otimes \wedge^{\bullet}\mathbf{C}^{5}$,
the tensor product of the defining representation of 
$\mathrm{SL}(\mathbf{C}^{2})$ and the exterior algebra 
of the defining representation of $\mathrm{SU}(5)$.
Under $\mathrm{SL}(\mathbf{C}^{2}) \times S(U(3) \times U(2))$,
this decomposes into 12 irreps corresponding to left and right handed 
electrons, neutrinos, up and down quarks and their antiparticles.

Unfortunately, $\mathrm{\textsl{diag}}\,(1,1,1,-1,-1) \in G$ acts
by $+1$ on right-handed fermions,
whereas $-\one \in \mathrm{SL}(\mathbf{C}^{2})$ acts by $-1$.
This means that $V$ does not define a representation of
$(\mathrm{SL}(\mathbf{C}^{2}) \times S(U(3) \times U(2))) / Z$
if one were to identify the central order 2 elements on both sides.

As the gauge group alone is of no use when trying to find
a $\spg$-structure,  
one has to involve the global $U(1)_{B} \times U(1)_{L}$-symmetries
connected to baryon and lepton number. (These rotate quarks and leptons 
independently.)

We conclude that the only infinitesimally natural $\spg$-structures 
relevant to the standard 
model are the ones associated to a homomorphism
\begin{equation} \label{homomo}
\pi_{1}(F^{+}(M)) \rightarrow 
S(U(3) \times U(2)) \times U(1)_{B} \times U(1)_{L}
\end{equation}
preserving $Z$, the subgroup of 
$U(1)_{B} \times U(1)_{L}$ generated by
$(-1,-1)$.

\subsubsection{Infinitesimally natural $\spg$-structures for 
the standard model}

Any manifold which possesses an infinitesimally natural 
$\spg$-structure for the standard model 
automatically permits an infinitesimally natural
$\mathrm{Spin}^{c}$-structure.
On the other hand, there do exist $\spg$-structures for the
standard model which are not $\mathrm{Spin}^{c}$.
We construct an example.

Consider de Sitter space 
$H = \{\vec{x} \in \mathbf{R}^{5} | 
- x^{2}_{0} + x^{2}_{1} + x^{2}_{2} + x^{2}_{3} + x^{2}_{4} = 1 \}$,
which has a pseudo-Riemannian metric $g$ with constant curvature 
induced by the Minkowski metric 
in the ambient $\mathbf{R}^{5}$.
Its group of orientation preserving isometries is $\mathrm{SO}(1,4)$,
and $H \simeq \mathrm{SO}(1,4)/ \mathrm{SO}(1,3)$.
Denote by $OF^{+\uparrow}_{g}(H)$ the bundle of orthogonal
frames with positive orientation and time-orientation.
By viewing $OF^{+\uparrow}_{g}(H)$ as a submanifold of 
$\mathbf{R}^{5} \times \mathrm{SO}(1,4)^{0}$, one can see that
$\mathrm{SO}(1,4)^{0}$ acts freely and transitively 
by $x : f \mapsto x_{*}f$. Therefore $OF^{+\uparrow}_{g}(H)$
is diffeomorphic to $\mathrm{SO}(1,4)^{0}$.

Now let $\Gamma < \mathrm{SO}(4)$ 
be a discrete group which acts freely, isometrically 
and properly discontinuously on $S^{3}$.
The manifold $\Gamma \backslash S^{3}$ is called 
a spherical space form. (See \cite{Wo} for a classification.)
As $\Gamma$ includes into $\mathrm{SO}(1,4)^{0}$, it acts
on $H$, making $M = \Gamma \backslash H$ into a pseudo-Riemannian manifold 
with constant curvature.

As $H$ is simply connected, we immediately see that $\pi_{1}(M) = \Gamma$.
Let us calculate the homotopy group of the frame bundle.
Because $OF_{g}^{+\uparrow}(M)$ is just 
$\Gamma \backslash OF_{g}^{+\uparrow}(H)$,
it equals
$\Gamma \backslash \mathrm{SO}(1,4)^{0}$. 
Going to the universal cover, we see that
$OF_{g}^{+}(M) = \tilde{\Gamma} \backslash \widetilde{\mathrm{SO}}(1,4)^{0}$.
As $\Gamma < \mathrm{SO}(4)$, we may consider $\tilde{\Gamma}$ to be 
the preimage of $\Gamma$ in $\mathrm{Spin}(4)$.
As the universal cover is simply connected, it is now clear that
$\pi_{1}(OF_{g}^{+\uparrow}(M)) = \tilde{\Gamma}$.
We get for free a homomorphism
$\tilde{\Gamma} \rightarrow \mathrm{Spin}(4)
\simeq
\mathrm{SU}(2)_{l} \times \mathrm{SU}(2)_{r}$, which maps the 
noncontractible loop in the fibre to $(-\one, -\one)$.

Triggered by the WMAP-data on cosmic background radiation, 
there has been some interest in the case where $\Gamma$ is $I^{*}$,
the binary icosahedral group \cite{Lea}. 
Some evidence against \cite{Kea} as well as in favour of \cite{Rea}
this hypothesis
appears to have been found. Without choosing sides in the debate, 
we simply point out that $I^{*} \backslash H$
gives rise to an interesting infinitesimally natural $\spg$-structure.

Under the identification
$\mathrm{Spin}(4) \simeq \mathrm{SU}(2)_{l} \times \mathrm{SU}(2)_{r}$,
we see that $\Gamma = I^{*} \times \one$ lives only in $\mathrm{SU}(2)_{l}$,
so that $\tilde{\Gamma}$ is the direct product of $I^{*} \times \one$
and the $\mathbf{Z}/2\mathbf{Z}$ generated by $(-\one, -\one)$.
One can therefore define a homomorphism (\ref{homomo})
by identifying $\mathrm{SU}(2)_{l}$ with 
$\mathrm{SU}(2)_{L}< G$, and mapping $(-\one, -\one)$ to 
$(-\one, -\one) \in U(1)_{B} \times U(1)_{L}$.
This yields an infinitesimally natural $\spg$-structure which uses the 
noncommutativity of the gauge group in an essential fashion. 
Note however that `ordinary'
$\mathrm{Spin}^{c}$-structures also exist on $M$.

\subsubsection{Extensions of the standard model}
The fact that $S(U(3) \times U(2))$ 
does not contribute to the obstruction of finding $\spg$-structures
on $M$
is due to the fact that
it never acts by $-1$ on $V$. This is not true for some GUT-type extensions of 
the standard model, such as the Pati-Salam
$\mathrm{SU}(2)_{L} \times \mathrm{SU}(2)_{R} \times \mathrm{SU}(4)$
model and anything which extends it, for example $\mathrm{Spin}(10)$.

If $N$ is the group of order 2 generated by $(-\one,-\one,-\one)$, then
infinitesimally natural $\spg$-structures in the Pati-Salam model 
correspond, neglecting global symmetries, to homomorphisms 
$\pi_{1}(F^{+}(M)) \rightarrow 
\mathrm{SU}(2)_{L} \times \mathrm{SU}(2)_{R} \times \mathrm{SU}(4)/N$
which take $Z$ to $\langle (-\one,-\one,\one)\rangle$.
This has the rather intriguing consequence that there may well
exist space-time manifolds $M$ which are compatible with the Pati-Salam model,
but not with the standard model.
A manifold $M$ would have this property
if 
the smallest quotient of $\pi_{1}(F(M))$ containing $Z$
is a nonabelian subgroup of 
$\mathrm{SU}(2)_{L} \times \mathrm{SU}(2)_{R} \times \mathrm{SU}(4)$
containing $(-\one,-\one,\one)$.

\section{Discussion}\label{acht}
A natural bundle is one in which diffeomorphisms of the base 
lift to automorphisms of the bundle in a local fashion.
In corollary \ref{PalaisTerng}, we have re\-derived the well known 
\cite{PT} 
result that any 
natural principal fibre bundle is associated to the $k^{\mathrm{th}}$
order frame bundle
$F^{k}(M)$. 

Bundles associated to spin structures are almost never natural, 
but nonetheless indispensable to physics.
This means that the notion of naturality will have to be extended.
One way to do this is to introduce
`gauge natural bundles' \cite{KMS}, 
which do not transform according to diffeomorphisms of the base $M$,
but according to automorphisms of a bundle over $M$ which has to be
specified.

We here propose to accommodate spin structures in a different fashion.
Although the Lorentz group does
not act on a spin structure, its Lie algebra does.
Analogously, we define an
`infinitesimally natural' bundle over $M$ to be a bundle
in which it is not the group of diffeomorphisms of $M$
that lifts, but its Lie algebra of vector fields.
%
This has the advantage that the link between the base and the 
fibres is not lost, so that a stress-energy-momentum tensor can be
constructed from Noether's principle \cite{No}. 

We assume only that the lift is a homomorphism of Lie algebras.
However,  
the careful analysis by Gotay and Marsden \cite{GM} 
reveals that in order to define a SEM-tensor from Noether's principle, 
one needs a lift that is given by a differential operator.
Assuming only that the lift is a homomorphism,
we prove that it is a differential operator in proposition \ref{alhetvet}.
This shows that 
one can construct a SEM-tensor from Noether's 
principle precisely when the bundle is infinitesimally natural.  
We therefore propose to describe fields on space-time 
only by sections of infinitesimally natural bundles. 

We have classified the infinitesimally natural principal fibre bundles.
According to 
theorem \ref{grondsmakreefknol},  
they are associated to the universal cover $\tf$ 
of the $k^{\mathrm{th}}$
order frame bundle, and their transformation behaviour is governed by
the disconnected group $G(k,M)$. 
This group depends on the manifold, and resembles a
$\mathrm{Pin}$ group
in that it can regulate parity and time reversal. 

The consistent description of fermions in the presence of a 
compact gauge group $G$ requires 
a $\spg$-structure rather than a spin structure \cite{HP}, \cite{AI}. 
All spin structures are infinitesimally natural, but some
$\spg$-structures are not. 
As described in theorem \ref{naturalspin},  
infinitesimally natural $\spg$-structures 
correspond to homomorphisms $\pi_{1}(F^{+}(M)) \rightarrow G$
that 
are injective on $\pi_{1}(\glnp)$.
 

The existence of infinitesimally natural $\spg$-structures therefore
provides 
an obstruction on space-time $M$ in terms of the symmetry group
$G$.
Some manifolds are even excluded for any (compact) choice of $G$. 
For example, $\mathbf{C}P^{2}$ does not admit
any infinitesimally natural $\spg$-structure 
because its frame bundle is simply connected. 
In our eyes, this disqualifies it as a model for space-time.


If all global symmetries are gauged, 
then theorem \ref{naturalspin} constitutes the
``connexion between the topology of space-time and the
spectrum of elementary particles'' alluded to in \cite{HP}.

\section*{Acknowledgements}
It is my pleasure to thank Marius Crainic, Hans Duistermaat, Johan van de Leur 
and Eduard Looijenga for their kind advice.

\end{document}